\documentclass{article}
\usepackage{amsmath,amsfonts,amssymb,color}
\usepackage{amsmath,amssymb,amstext,latexsym,amsfonts,graphicx,fancybox}
\newtheorem{lem}{Lemma}[section]

\newtheorem{rem}[lem]{Remark}
\newtheorem{theo}[lem]{Theorem}

\newcommand{\N}{\mathbb{N}}
\newcommand{\R}{\mathbb{R}}

\newcommand{\E}{\mathbb{E}}

\newcommand{\dL}{\mathbb{L}}
\newcommand{\C}{\mathbb{C}}

\newcommand{\p}{\mathbb{P}}

\newcommand{\pt}{\!\!\!\bf .}

\newcommand{\ind}{\mbox{1}\kern-.21em \mbox{I}}

\newcommand{\vp}{\varphi}
\newcommand{\cN}{\mathcal{N}}

\newcommand{\cD}{\mathcal{D}}

\newcommand{\cL}{\mathcal{L}}
\newcommand{\cH}{\mathcal{H}}

\newcommand{\cR}{\mathcal{R}}

\newcommand{\cO}{\mathcal{O}}
\font\calcal=cmsy10 scaled\magstep1
\def\build#1_#2^#3{\mathrel{\mathop{\kern 0pt#1}\limits_{#2}^{#3}}}
\def\liml{\build{\longrightarrow}_{}^{{\mbox{\calcal L}}}}

\def\videbox{\mathbin{\vbox{\hrule\hbox{\vrule height1ex \kern.5em\vrule height1ex}\hrule}}}

\def\es#1{{\mathbb E}\left[{#1}\right]}

\def\loga#1{\log \left[ {#1} \right]}

\def\demend{\hfill $\videbox$\\}

\DeclareMathOperator{\re}{Re}
\DeclareMathOperator{\im}{Im}
\DeclareMathOperator{\Arg}{Arg}

\begin{document}
\title{Sharp large deviations for the non-stationary Ornstein-Uhlenbeck process}
\date{}
\author{Bernard Bercu
\thanks{Universit\'e Bordeaux 1, Institut de Math\'ematiques de Bordeaux,
UMR C5251, 351 cours de la lib\'eration, 33405 Talence cedex, France. e-mail:
Bernard.Bercu@math.u-bordeaux1.fr}
\and Laure Coutin 
\thanks{Universit\'e Paul Sabatier, Institut de Math\'ematiques de Toulouse,
UMR C5583, 31062 Toulouse Cedex 09, France. e-mail:
laure.coutin@math.univ-toulouse.fr}
\and Nicolas Savy
\thanks{Universit\'e Paul Sabatier, Institut de Math\'ematiques de Toulouse,
UMR C5583, 31062 Toulouse Cedex 09, France. e-mail:
nicolas.savy@math.univ-toulouse.fr}}
\maketitle
\[\hbox{{\bf Abstract}}\]
For the Ornstein-Uhlenbeck process, the asymptotic behavior of the maximum likelihood estimator of the drift parameter
is totally different in the stable, unstable, and explosive cases. Notwithstanding of this trichotomy, we investigate sharp 
large deviation principles for this estimator in the three situations. 
In the explosive case, we exhibit a very unusual rate function with a shaped flat valley 
and an abrupt discontinuity point at its minimum.
\vspace{2ex} \newline
{\bf A.M.S. Classification:} 60F10, 60G15, 62A10.
\newline
{\bf Key words:} Large deviations, Ornstein-Uhlenbeck process, Likelihood estimation.

\section{Introduction.}
\renewcommand{\theequation}{\thesection.\arabic{equation}}
\setcounter{section}{1}
\setcounter{equation}{0}

Consider the Ornstein-Uhlenbeck process observed over the time interval $[0,T]$
\begin{equation}
\label{OUPROCESS}
dX_t = \theta X_t dt + dB_t
\end{equation}
where $(B_t)$ is a standard Brownian motion and the drift $\theta$ is an unknown real parameter. For the sake of
simplicity, we choose the initial state $X_0 = 0.$ The process is said to be stable if $\theta < 0$, unstable if $\theta = 0$, and explosive if $\theta > 0$. 
The maximum likelihood estimator of $\theta$ is given by
\begin{equation} 
\label{DEFTHETAHAT}
\widehat{\theta}_T = \frac{\int_0^T X_t dX_t}{\int_0^T X_t^2 dt} = \frac{X_T^2 - T}{2\int_0^T X_t^2 dt}.
\end{equation}
It is well-known  (see e.g. \cite{LS01} page 234) that in the stable, unstable, and explosive cases
\begin{equation*} 
\lim_{T\rightarrow \infty}\widehat{\theta}_T = \theta \hspace{1cm}\text{a.s.}
\end{equation*}
However, the asymptotic normality is totally different in the three situations. As a matter of fact, if $\theta <0$, the process
$(X_T)$ is positive recurrent and Brown and Hewitt \cite{BH75} have shown the asymptotic normality
\begin{equation*}
\sqrt{T}(\widehat{\theta}_T - \theta) \liml \cN(0, -2 \theta).
\end{equation*}
Moreover, if $\theta =0$, the process $(X_T)$ is null recurrent and it was proved by Feigin \cite{F79} that
\begin{equation*}
T(\widehat{\theta}_T -\theta) \liml \frac{\int_0^1 W_t dW_t}{\int_0^1 W_t^2 dt} = \frac{W_1^2 - 1}{2\int_0^1 W_t^2 dt}
\end{equation*}
where $(W_t)$ is a standard Brownian motion. Furthermore,  if $\theta >0$, the process $(X_T)$ is transient and we know from Feigin \cite{F76},
(see also \cite{K04} page 304), that
\begin{equation*}
\exp(\theta T)(\widehat{\theta}_T -\theta)  \liml 2 \theta \, \Bigl(\frac{Y}{Z}\Bigr)
\end{equation*}
where $Y, Z$ are two independent Gaussian $\cN(0,1)$ random variables which implies that the limiting ratio $Y/Z$
has a Cauchy distribution. Notwithstanding of this trichotomy, our goal is to establish the large
deviation properties for $(\widehat{\theta}_T)$ in the stable, unstable, and explosives cases.
We refer the reader to the excellent book by Dembo and Zeitouni \cite{DZ} on the theory of large deviations.
First of all, in the stable case, Florens-Landais and Pham \cite{FP99} proved the following large deviation principle
(LDP) for $(\widehat{\theta}_T)$.

\begin{lem}$\pt$ 
\label{L-STA}
If $\theta <0$, then $(\widehat{\theta}_T)$ satisfies an LDP with speed $T$ and good rate function
\begin{eqnarray}
\label{RATESTA}
I(c)= \left \{ \begin{array}{lll}
    {\displaystyle \hspace{-1ex} -\frac{(c-\theta)^2}{4c}} & \text{ if } \hspace{2ex} {\displaystyle c<\frac{\theta}{3}}, \vspace{1ex}\\
     2c-\theta & \text{ otherwise. }     \end{array}  \right.
\end{eqnarray}
\end{lem}

This result was extended by the  following sharp large deviation principle
(SLDP) for $(\widehat{\theta}_T)$ established  by Bercu and Rouault \cite{BR01}. 

\begin{theo}$\pt$
\label{T-SLDPSTA}
Consider the Ornstein-Uhlenbeck process given by (\ref{OUPROCESS}) where the drift parameter $\theta <0$. 
\newline
a) For all $c < \theta$, there exists a sequence $(b_{c,k})$ such that, for any $p> 0$ and $T$ large enough
\begin{equation}
\label{SLDPSTA1}
\p(\widehat{\theta}_T \leq c) = \frac{-\exp(-TI(c) + H(a_c))}{a_c \sigma_c \sqrt{2 \pi T}}
\left[1 + \sum_{k=1}^{p}\frac{b_{c,k}}{T^k} + {\cal O}\Bigl(\frac{1}{T^{p+1}}\Bigr)\right]
\end{equation}
where
\begin{equation}
\label{DEFAC1}
a_c=\frac{c^2-\theta^2}{2c}\hspace{1cm}\text{and} \hspace{1cm} \sigma_c^2=-\frac{1}{2c}
\end{equation}
\begin{equation}
H(a_c)=-\frac{1}{2} \log \left(\frac{(c+\theta)(3c-\theta)}{4c^2}\right)
\label{DEFHAC1}
\end{equation}
while, for all $\theta < c < \theta/3$,
\begin{equation}
\label{SLDPSTA2}
\p(\widehat{\theta}_T \geq c) = \frac{\exp(-TI(c) + H(a_c))}{a_c \sigma_c \sqrt{2 \pi T}}
\left[1 + \sum_{k=1}^{p}\frac{b_{c,k}}{T^k} + {\cal O}\Bigl(\frac{1}{T^{p+1}}\Bigr)\right].
\end{equation}
\vspace{-2ex}
\newline
b) For all $c>\theta/3$ with $c \neq 0$, there exists a 
sequence $(d_{c,k})$ such that, for any $p>0$ and $T$ large enough
\begin{equation}
\label{SLDPSTA3}
\p(\widehat{\theta}_T \geq c) = 
\frac{\exp(-TI(c) + K(c))}{a_c \sigma_c\sqrt{2 \pi T}}
\left[1 + \sum_{k=1}^{p}\frac{d_{c,k}}{T^k}
+ {\cal O}\Bigl(\frac{1}{T^{p+1}}\Bigr)\right]
\end{equation}
where 
\begin{equation}
\label{DEFAC2}
a_c=2(c-\theta)\hspace{1cm}\text{and} \hspace{1cm}\sigma_c^2= \frac{c^2}{2(2c-\theta)^3}
\end{equation}
\begin{equation}
K(c)=-\frac{1}{2} \log \left(\frac{(c-\theta)(3c-\theta)}{4c^2}\right).
\label{DEFKC2}
\end{equation}
\vspace{-2ex}
\newline
c) For $c=\theta/3$, there exists a sequence $(e_{k})$
such that, for any $p>0$ and $T$ large enough
\begin{equation}
\label{SLDPSTA4}
\p(\widehat{\theta}_T \geq c) = 
\frac{÷exp(-TI(c))}{2 \pi T^{1/4}}
\frac{\Gamma(1/4)}{a_{\theta}^{3/4}\sigma_{\theta}}
\left[1 + \sum_{k=1}^{2p}\frac{e_{k}}{(\sqrt{T})^k}
+ {\cal O}\Bigl(\frac{1}{T^{p}\sqrt{T}}\Bigr)\right]
\end{equation}
where
\begin{equation}
a_{\theta}= -\frac{4\theta}{3}
\hspace{1cm}\text{and} \hspace{1cm}
\sigma_{\theta}^2=-\frac{3}{2\theta}.
\end{equation}
\vspace{-2ex}
\newline
d) Finally, for $c=0$, $p>0$ and for $T$ large enough
\begin{equation}
\label{SLDPSTA5}
\p(\widehat{\theta}_T \geq 0) = 2
\frac{\exp(-TI(c))}{\sqrt{2 \pi T}\sqrt{-2\theta}}
\left[1 + \sum_{k=1}^{p}\frac{(2k)!}{2^{2k}\theta^kT^k k!}
+ {\cal O}\Bigl(\frac{1}{T^{p+1}}\Bigr)\right].
\end{equation}
\end{theo}

Our purpose is to extend this investigation by establishing SLDP for $(\widehat{\theta}_T)$ in the explosive and unstable cases.
Similar results in discrete time for the Gaussian autoregressive process may be found in \cite{B01}.
We also refer the reader to \cite{BCS11} where SLDP for the maximum likelihood estimator
of $\theta$ is proved for the stable Ornstein-Uhlenbeck process driven by a fractional Brownian motion.
We wish to mention here that it could be possible to extend the previous work of Zani \cite{Z02} to
generalized squared radial Ornstein-Uhlenbeck processes with parameter $\theta >0$.

\section{A keystone lemma.}
\setcounter{section}{2}
\setcounter{equation}{0}


The sharp large deviations properties of $(\widehat{\theta}_T)$ are closely related with the ones of
$$
Z_T(c) = \int_0^T X_t dX_t - c \int_0^T X_t^2 dt
$$
with $c \in \R$ since $\p(\widehat{\theta}_T \geq c) = \p(Z_T(c) \geq 0)$. One has to keep in mind that the threshold 
$c$ for $\widehat{\theta}_T$ appears like a parameter for $Z_T$. Denote by $\cL_T$ 
the normalized cumulant generating function of 
$Z_T(c)$
$$
\cL_T(a) = \frac{1}{T} \log \E\Bigl[\exp(a Z_T(c))\Bigr]
$$
where the parameter $c$ is omitted in order to simplify the notation. All our analysis relies on the following keystone lemma which is true
as soon as the drift parameter $\theta \geq 0$.

\begin{lem} $\pt$\label{L-LEMFOND}
Let $\Delta_c = \{ a \in \R,  ~\theta^2 +2ac > 0,~a + \theta < \sqrt{\theta^2+2ac} \}$ be the effective domain of the limit $\cL$ of $\cL_T$ and set $\vp(a) = - \sqrt{\theta^2+2ac}$, $\tau(a) = a + \theta - \vp(a)$ and 
$h(a) = (a + \theta)/\vp(a)$.
\newline
a) For all $a \in \Delta_c$, we have 
\begin{equation} 
\label{E-DECFOND}
\cL_T(a)=\cL(a) + \frac{1}{T} \cH(a) + \frac{1}{T} \cR_T(a)
\end{equation}
where
\begin{eqnarray}
\cL(a)   &=& -\frac{1}{2} \left(a + \theta + \sqrt{\theta^2+2ac} \right),\label{E-DEFL}\\
\cH(a)   &=& -\frac{1}{2} \log \left( \frac{1}{2} (1+ h(a)) \right),\label{E-DEFH}\\
\cR_T(a) &=&-\frac{1}{2} \log \left( 1 + \frac{1 - h(a)}{ 1 + h(a)} \exp(2 \varphi(a) T) \right) \label{E-DEFR}.
\end{eqnarray}
\newline
b) Moreover, the remainder $\cR_T(a)$ goes to zero exponentially fast as
\begin{equation} 
\label{CONTRT}
\cR_T(a) = \cO \Bigl( \exp(2 \varphi(a) T) \Bigr).
\end{equation}
\end{lem}

\noindent{\bf Proof.} The proof of Lemma \ref{L-LEMFOND} is given in Appendix A.\demend

\section{Sharp Large deviations results.}
\setcounter{section}{3}
\setcounter{equation}{0}


We shall now focus our attention on the explosive case $\theta >0$. It immediately follows from 
(\ref{OUPROCESS}) that 
\begin{equation}
\label{SOLOUP}
X_T=\exp(\theta T) \int_0^T \exp(-\theta t)dB_t
\end{equation}
leading to $\exp(-\theta T)X_T= Y_T$ where
\begin{equation*}
Y_T=\int_0^T \exp(-\theta t)dB_t.
\end{equation*}
The Gaussian process $(Y_T)$ converges almost surely and in mean square to the Gaussian nondegenerate random variable
\begin{equation*}
Y=\int_0^\infty \exp(-\theta t)dB_t.
\end{equation*}
Hence, via Toeplitz's lemma
\begin{equation*}
\lim_{T\rightarrow \infty} \frac{1}{\exp(2 \theta T)} \int_0^T X_t^2dt= \frac{Y^2}{2\theta} \hspace{1cm}\text{a.s.}
\end{equation*}
Consequently, one can expect for $(\widehat{\theta}_T)$ an LDP with speed $\exp(2 \theta T)$. However, 
$(\widehat{\theta}_T)$ is a sequence of self-normalized random variables and we shall show that
$(\widehat{\theta}_T)$ satisfies an LDP similar to that of Lemma \ref{L-STA} with speed $T$.

\begin{lem}$\pt$ 
\label{L-LDPE}
If $\theta >0$, then $(\widehat{\theta}_T)$ satisfies an LDP with speed $T$ and good rate function
\begin{eqnarray}
\label{RATE}
I(c)= \left \{ \begin{array}{lll}
    {\displaystyle \hspace{-1ex} -\frac{(c - \theta)^2}{4c} } & \text{ if } \hspace{2ex}  c \leq -\theta, \vspace{1ex}\\
    \theta  & \text{ if } \hspace{2ex}|c| < \theta, \vspace{2ex}\\
     0   & \text{ if } \hspace{2ex} c = \theta, \vspace{2ex}\\
     2c - \theta  &\text{ if } \hspace{2ex} c > \theta.
\end{array}  \right.
\end{eqnarray}
\end{lem}

\begin{rem}
As for the Gaussian autoregressive process \cite{B01}, one can observe that the rate function $I$ in the explosive case
is really unusual with a shaped flat valley and an abrupt discontinuity point at its minimum. 
It is possible to give some intuition on the size of the discontinuity jump. As a matter of fact, 
we already saw in the introduction that, if $\theta >0$, 
\begin{equation*}
\exp(\theta T)(\widehat{\theta}_T -\theta)  \liml 2 \theta \, \Bigl(\frac{Y}{Z}\Bigr)
\end{equation*}
where $Y, Z$ are two independent Gaussian $\cN(0,1)$ random variables. The size of the jump is precisely given
by the logarithm of the rate $\exp(\theta T)$ properly normalized, 
\begin{equation*}
\frac{1}{T}\log(\exp(\theta T))=\theta.
\end{equation*}
\end{rem}

The SLDP for $(\widehat{\theta}_T)$, quite similar to the one established in the stable case,  is as follows. 

\begin{theo}$\pt$
\label{T-SLDPE}
Consider the Ornstein-Uhlenbeck process given by (\ref{OUPROCESS}) where the drift parameter $\theta >0$. 
\newline
a) For all $c < -\theta$, there exists a sequence $(b_{c,k})$ such that, for any $p> 0$ and $T$ large enough
\begin{equation}
\label{SLDPE1}
\p(\widehat{\theta}_T \leq c) = \frac{-\exp(-TI(c) + H(a_c))}{a_c \sigma_c \sqrt{2 \pi T}}
\left[1 + \sum_{k=1}^{p}\frac{b_{c,k}}{T^k} + {\cal O}\Bigl(\frac{1}{T^{p+1}}\Bigr)\right]
\end{equation}
where
\begin{equation}
\label{DEFACE1}
a_c=\frac{c^2-\theta^2}{2c}\hspace{1cm}\text{and} \hspace{1cm} \sigma_c^2=-\frac{1}{2c}
\end{equation}
\begin{equation}
H(a_c)=-\frac{1}{2} \log \left(\frac{(c+\theta)(3c-\theta)}{4c^2}\right).
\label{DEFHACE1}
\end{equation}
\vspace{-2ex}
\newline
b) For all $c>\theta$, there exists a sequence $(d_{c,k})$ such that, for any $p>0$ and $T$ large enough
\begin{equation}
\label{SLDPE2}
\p(\widehat{\theta}_T \geq c) = 
\frac{\exp(-TI(c) + K(c))}{a_c \sigma_c\sqrt{2 \pi T}}
\left[1 + \sum_{k=1}^{p}\frac{d_{c,k}}{T^k}
+ {\cal O}\Bigl(\frac{1}{T^{p+1}}\Bigr)\right]
\end{equation}
where 
\begin{equation}
\label{DEFACE2}
a_c=2(c-\theta)\hspace{1cm}\text{and} \hspace{1cm}\sigma_c^2= \frac{c^2}{2(2c-\theta)^3}
\end{equation}
\begin{equation}
K(c)=-\frac{1}{2} \log \left(\frac{(c-\theta)(3c-\theta)}{4c^2}\right).
\label{DEFKCE2}
\end{equation}
\vspace{-2ex}
\newline
c) For all $|c|<\theta$  with $c \neq 0$, there exists a sequence $(e_{c,k})$
such that, for any $p>0$ and $T$ large enough
\begin{equation}
\label{SLDPE3}
\p(\widehat{\theta}_T \leq c) = 
\frac{\exp(-TI(c) + J(c))}{a_c \sigma_c\sqrt{2 \pi T}}
\left[1 + \sum_{k=1}^{p}\frac{e_{c,k}}{T^k}
+ {\cal O}\Bigl(\frac{1}{T^{p+1}}\Bigr)\right]
\end{equation}
where 
\begin{equation}
\label{DEFACE3}
a_c=\frac{\theta}{c+\theta}\hspace{1cm}\text{and} \hspace{1cm}\sigma_c^2= \frac{c^2}{2\theta^3}
\end{equation}
\begin{equation}
J(c)=-\frac{1}{2} \log \left(\frac{(\theta-c)(\theta+c)}{4c^2}\right).
\label{DEFJCE3}
\end{equation}
\vspace{-2ex}
\newline
d) For $c=-\theta$, there exists a sequence $(f_{k})$
such that, for any $p>0$ and $T$ large enough
\begin{equation}
\label{SLDPE4}
\p(\widehat{\theta}_T \leq c) = 
\frac{\exp(-TI(c))}{2 \pi T^{1/4}}
\frac{\Gamma(1/4)}{a_{\theta}^{3/4}\sigma_{\theta}}
\left[1 + \sum_{k=1}^{2p}\frac{f_{k}}{(\sqrt{T})^k}
+ {\cal O}\Bigl(\frac{1}{T^{p}\sqrt{T}}\Bigr)\right]
\end{equation}
where
\begin{equation}
\label{DEFAC4}
a_{\theta}= \sqrt{\theta}
\hspace{1cm}\text{and} \hspace{1cm}
\sigma_{\theta}^2=\frac{1}{2\theta}.
\end{equation}
\vspace{-2ex}
\newline
e) Finally, for $c=0$, $p>0$ and for $T$ large enough
\begin{equation}
\label{SLDPE5}
\p(\widehat{\theta}_T \leq 0) = 2
\frac{\exp(-T I(c)) \sqrt{2 \theta T}}{\sqrt{2 \pi}}
\left[1 + \sum_{k=1}^{p}\frac{(-1)^{k} (\theta T e^{-2\theta T})^k}{(2k+1) k!} 
+ {\cal O}\Bigl((Te^{-2\theta T})^{p+1}\Bigr)\right].
\end{equation}
\end{theo}

\begin{rem} $\pt$
One can observe that all the sequences $(b_{c,k})$ $(d_{c,k})$, $(e_{c,k})$ may be explicitly calculated as in Theorem 4.1 of \cite{BR01}.
\end{rem}

\noindent{\bf Proof.} The proofs are given in Section \ref{S-PROOF}. \demend

The unstable case $\theta = 0$ can be handled exactly as the explosive case $\theta > 0$ since Lemma \ref{L-LEMFOND}
is also true in the unstable situation. Consequently, we directly obtain the LDP and SLDP for $(\widehat{\theta}_T)$ in the unstable case
by replacing $\theta$ by $0$ in the previous results.

\begin{lem}$\pt$ 
\label{L-LDPUNS}
If $\theta =0$, then $(\widehat{\theta}_T)$ satisfies an LDP with speed $T$ and good rate function
\begin{eqnarray}
\label{RATEUNS}
I(c)= \left \{ \begin{array}{lll}
    {\displaystyle \hspace{-1ex} -\frac{c}{4} } & \text{ if } \hspace{2ex}  c \leq 0, \vspace{1ex}\\
     2c   &\text{ otherwise}.
\end{array}  \right.
\end{eqnarray}
\end{lem}

\begin{theo}$\pt$
\label{T-SLDPUNS}
Consider the Ornstein-Uhlenbeck process given by (\ref{OUPROCESS}) where the drift parameter $\theta =0$. 
\newline
a) For all $c < 0$, there exists a sequence $(b_{c,k})$ such that, for any $p> 0$ and $T$ large enough
\begin{equation}
\label{SLDPUNS1}
\p(\widehat{\theta}_T \leq c) = \frac{-2\exp(-TI(c))}{a_c \sigma_c \sqrt{6 \pi T}}
\left[1 + \sum_{k=1}^{p}\frac{b_{c,k}}{T^k} + {\cal O}\Bigl(\frac{1}{T^{p+1}}\Bigr)\right]
\end{equation}
where $a_c=c/2$ and $\sigma_c^2=-1/(2c)$.
\newline
b) For all $c>0$, there exists a sequence $(d_{c,k})$ such that, for any $p>0$ and $T$ large enough
\begin{equation}
\label{SLDPUNS2}
\p(\widehat{\theta}_T \geq c) = 
\frac{2\exp(-TI(c))}{a_c \sigma_c\sqrt{6 \pi T}}
\left[1 + \sum_{k=1}^{p}\frac{d_{c,k}}{T^k}
+ {\cal O}\Bigl(\frac{1}{T^{p+1}}\Bigr)\right]
\end{equation}
where $a_c=2c$ and $\sigma_c^2=1/(16c)$.
\end{theo}

\section{Proofs of the main results.} \label{S-PROOF}
\setcounter{section}{4}
\setcounter{equation}{0}


\subsection{Proof of Theorem \ref{T-SLDPE} a).}

We first focus our attention on the easy case $c < - \theta$. One can observe that $a_c$, given by \eqref{DEFACE1}, 
belongs to the effective domain $\Delta_c=]-\infty,0[$ whenever $c < - \theta$. Consider the usual change of probability
\begin{equation} \label{PT1} 
\frac{d\p_T}{d\p} = \exp\Bigl( a_c Z_T(c) - T \cL_T(a_c) \Bigr)
\end{equation}
and denote by $\E_T$ the expectation associated with $\p_T$. We clearly have
\begin{eqnarray*}
\p(\widehat{\theta}_T \leq c) 
&=& \p(Z_T(c) \leq 0)= \E[\ind_{Z_T(c) \leq 0}],\\
&=& \E_T \Bigl[\exp( - a_c Z_T(c) + T \cL_T(a_c) ) \ind_{Z_T(c) \leq 0}\Bigr],\\
&=& \exp\Bigl( T \cL_T(a_c) \Bigr) \E_T \Bigl[ \exp( - a_c Z_T(c) ) \ind_{Z_T(c) \leq 0}\Bigr].
\end{eqnarray*}
Consequently, we can split $\p(\widehat{\theta}_T \leq c)$ into two terms, $\p(\widehat{\theta}_T \leq c)=A_T B_T$ with
\begin{eqnarray} 
A_T &=& \exp(T \cL_T(a_c)), \label{AT1}\\
B_T &=& \E_T[\exp(-a_c Z_T(c))\ind_{Z_T(c) \leq 0}]. \label{BT1}
\end{eqnarray}
On the one hand, we can deduce from \eqref{E-DECFOND} and \eqref{CONTRT} together with the definition \eqref{RATE}  of $I$ that
\begin{eqnarray}
A_T  &=& \exp\Bigl(T \cL(a_c) + \cH(a_c)+\cR_T(a_c)\Bigr), \nonumber \\
\label{DEVAT1}
A_T  &=& \exp\Bigl(- T I(c) + \cH(a_c)\Bigr) \left( 1 + \cO \left(e^{2 T c} \right) \right) .
\end{eqnarray}
It only remains to provide the expansion for $B_T$.

\begin{lem}$\pt$ \label{L-BT1}
For all $c < -\theta$, there exists a sequence $(\beta_k)$ such that, for any $p>0$ and $T$ large enough,
\begin{equation} \label{DEVBT1}
B_T= \frac{\beta_0}{\sqrt{T}}\left[ 1+\sum_{k=1}^p \frac{\beta_k}{T^k} + {\mathcal O}\Bigl( \frac{1}{T^{p+1}}\Bigr)\right].
\end{equation}
The sequence $(\beta_k)$ only depends on the derivatives of $\cL$ and $\cH$ evaluated at point $a_c$. For example,
\begin{equation*}
\beta_0=-\frac{1}{a_c\sigma_c\sqrt{2\pi}}.
\end{equation*}
\end{lem}

\noindent
{\bf Proof.} The proof of Lemma \ref{L-BT1} is given in Appendix C.\demend

\noindent
{\bf Proof of Theorem \ref{T-SLDPE} a).} The expansion \eqref{SLDPE1} immediately follows from \eqref{DEVAT1} and \eqref{DEVBT1}.
\demend

\subsection{Proof of Theorem \ref{T-SLDPE} b).}

In the more complicated case $c > \theta$, the effective domain $\Delta_c=]0,2(c-\theta)[$ and the function $\cL$ is decreasing
over the interval $]0,2(c-\theta)[$ as
\begin{equation*} 
\cL^{\prime}(a) = -\frac{1}{2} \left( 1 + \frac{c}{\sqrt{\theta^2 + 2ac}}\right).
\end{equation*}
Consequently, $\cL$ reaches its minimum at the value $a_c=2(c-\theta)$ given by \eqref{DEFACE2}.
Therefore, it is necessary to make use of a slight modification of the strategy of time varying change of probability 
proposed by Bryc and Dembo \cite{BD}. The key point is that there exists a unique $a_T$, which belongs to the interior of $\Delta_c$ and converges to its border $a_c = 2(c-\theta)$ as $T$ goes to infinity, solution of some suitable implicit equation. Hereafter, we introduce the new probability measure
\begin{equation} \label{PT2}
\frac{d \p_T}{d \p }=\exp \Bigl( a_T Z_T(c) - T \cL_T(a_T) \Bigr)
\end{equation}
and we denote by $\E_T$ the expectation under $\p_T$. It leads to the decomposition $\p(\widehat{\theta}_T \geq c)=A_TB_T$
where
\begin{eqnarray}
A_T &= & \exp \left( T \cL_T(a_T) \right),\label{AT2}\\
B_T &=& \E_T \Bigl[\exp ( - a_T Z_T(c) ) \ind_{Z_T(c) \geq 0}\Bigr].\label{BT2}
\end{eqnarray}
The proof now splits into two parts, the first one is devoted to the expansion of $A_T$ while the second one gives the expansion for $B_T$.\\

\begin{lem}$\pt$
\label{L-AT2}
For all $c > \theta$, there exists a unique $a_T$, which belongs to the interior of $\Delta_c$ and converges to its border $a_c=2(c-\theta)$
as $T$ goes to infinity, solution of the implicit equation
\begin{equation}
\label{IMPEQ}
\cL^{\prime}(a) + \frac{1}{T}\cH^{\prime}(a) = 0
\end{equation}
where the functions $\cL$ and $\cH$ are given by \eqref{E-DEFL} and \eqref{E-DEFH}.
Moreover, there exists a sequence $(\gamma_k)$ such that, for any $p > 0$ and $T$ large enough,
\begin{equation}
A_T = \exp \left(- T I(c) + P(c) \right) \sqrt{ eT } \left[ 1 + \sum_{k=1}^p \frac{\gamma_k}{T^k} + {\mathcal O}\Bigl( \frac{1}{T^{p+1}}\Bigr) \right] \label{DEVATFIN2}
\end{equation}
where
\begin{equation} \label{DEFPQ1}
P(c) = -\frac{1}{2} \log \left( \frac{(c - \theta)}{2(2c- \theta)(3c- \theta)}\right).
\end{equation}
The sequence $(\gamma_k)$ only depends on the Taylor expansion of $a_T$ at the neighborhood of $a_c$ together with the derivatives of $\cL$ 
and $\cH$ at point $a_c$. For example,
\begin{equation*}
\gamma_1 = \frac{c(c^2-3\theta c +\theta^2)}{2(c - \theta)( \theta - 2c)(3c-\theta)^2}.
\end{equation*}
\end{lem}

\noindent {\bf Proof.} The proof of Lemma \ref{L-AT2} is given in Appendix B.\demend

\noindent
It now remains to give the expansion for $B_T$.

\begin{lem}$\pt$ \label{L-BT2}
For all $c > \theta$, there exists a sequence $(\delta_k)$ such that, for any $p>0$ and $T$ large enough,
\begin{equation} \label{DEVBT2}
B_T= \sum_{k=1}^p \frac{\delta_k}{T^k} + {\mathcal O}\Bigl( \frac{1}{T^{p+1}}\Bigr).
\end{equation}
The sequence $(\delta_k)$ only depends on the Taylor expansion of $a_T$ at the neighborhood of $a_c$ together with the derivatives of $\cL$ and
$\cH$ at point $a_c$. For example,
\begin{equation*}
\delta_1=\frac{1}{a_c \delta \sqrt{2 \pi e}}
\hspace{1cm}\text{where} \hspace{1cm}
\delta=- \cL^{\prime}(a_c)=\frac{(3c-\theta)}{2(2c-\theta)}.
\end{equation*}
\end{lem}

\noindent
{\bf Proof.} The proof of Lemma \ref{L-BT2} is given in Appendix C.\demend

\noindent
{\bf Proof of Theorem \ref{T-SLDPE} b).} The expansions \eqref{DEVATFIN2} and \eqref{DEVBT2} immediately imply \eqref{SLDPE2}.
\demend

\subsection{Proof of Theorem \ref{T-SLDPE} c).}

In the case $ |c| < \theta$ and $c \ne 0$, one can easily see that the effective domain is
\begin{eqnarray*}
\Delta_c = \left \{ \begin{array}{lll}
 ]-\infty, 0[ & \text{ if }   -\theta < c < 0, \vspace{1ex} \\
   {\displaystyle  ]-\frac{\theta^2}{2c}, 0[ } & \text{ if } \hspace{2ex} {\displaystyle 0 < c \leq \frac{\theta}{2}}, \vspace{1ex} \\
     ]2(c - \theta), 0[   & \text{ if } \hspace{2ex} {\displaystyle \frac{\theta}{2} \leq c < \theta}. 
\end{array}  \right.
\end{eqnarray*}
In addition, the function $\cL$ is always decreasing over $\Delta_c$ and $\cL$ reaches its minimum at the origin.
Consequently, the proof follows essentially the same lines as the one for $c > \theta$ with $a_c = 0$. In fact, 
with the new probability measure given by \eqref{PT2}, we have the decomposition $\p(\widehat{\theta}_T \leq c)=A_TB_T$
where
\begin{eqnarray}
A_T &=& \exp \left( T \cL_T(a_T) \right),\label{AT3}\\
B_T &=& \E_T \Bigl[\exp ( - a_T Z_T(c) ) \ind_{Z_T(c) \leq 0}\Bigr].\label{BT3}
\end{eqnarray} 
The proof is also divided into two parts, the first one is devoted to the expansion of $A_T$ while the second one gives the expansion for $B_T$.\\
\begin{lem} \label{L-AT3}
For all $ |c| < \theta$ and $c \ne 0$, there exists a unique $a_T$, which belongs to the interior of $\Delta_c$ and converges to the origin 
as $T$ goes to infinity, solution of the implicit equation
\begin{equation}
\label{IMPEQN}
\cL^{\prime}(a) + \frac{1}{T}\cH^{\prime}(a) = 0
\end{equation}
where the functions $\cL$ and $\cH$ are given by \eqref{E-DEFL} and \eqref{E-DEFH}.
Moreover, there exists a sequence $(\gamma_k)$ such that, for any $p > 0$ and $T$ large enough,
\begin{equation}
A_T = \exp \left(- T I(c) + P(c) \right) \sqrt{ eT } \left[ 1 + \sum_{k=1}^p \frac{\gamma_k}{T^k} + {\mathcal O}\Bigl( \frac{1}{T^{p+1}}\Bigr) \right] \label{DEVATFIN3}
\end{equation}
where
\begin{equation} \label{DEFPQ2}
P(c) =  -\frac{1}{2} \log \left( \frac{(\theta-c)}{2\theta(c + \theta)}\right).
\end{equation}
The sequence $(\gamma_k)$ only depends on the Taylor expansion of $a_T$ at the neighborhood of the origin together with the derivatives of $\cL$ and
$\cH$ at $0$. For example,
\begin{equation*}
\gamma_1 = - \frac{c(c^2 + \theta c - \theta^2)}{2\theta(c - \theta) (c + \theta)^2} .
\end{equation*}
\end{lem}

\noindent {\bf Proof.} The proof of Lemma \ref{L-AT3} is given in Appendix B.\demend

\noindent
The expansion of $B_T$ is as follows.

\begin{lem}$\pt$ \label{L-BT3}
For all $ |c| < \theta$ and $c \ne 0$, there exists a sequence $(\delta_k)$ such that, for any $p>0$ and $T$ large enough,
\begin{equation} \label{DEVBT3}
B_T= \sum_{k=1}^p \frac{\delta_k}{T^k} + {\mathcal O}\Bigl( \frac{1}{T^{p+1}}\Bigr).
\end{equation}
The sequence $(\delta_k)$ only depends on the Taylor expansion of $a_T$ at the neighborhood of the origin together with the derivatives of $\cL$ and
$\cH$ at $0$. For example,
\begin{equation*}
\delta_1=\frac{1}{a_c \delta \sqrt{2 \pi e}}
\hspace{1cm}\text{where} \hspace{1cm}
\delta=- \cL^{\prime}(0)=-\frac{(c+\theta)}{2\theta}.
\end{equation*}
\end{lem}

\noindent
{\bf Proof.} The proof of Lemma \ref{L-BT3} is given in Appendix C.\demend

\noindent
{\bf Proof of Theorem \ref{T-SLDPE} c).} The expansion \eqref{SLDPE3} immediately follows from  \eqref{DEVATFIN3} and \eqref{DEVBT3}.
\demend

\subsection{Proof of Theorem \ref{T-SLDPE} d). }

In the particular case $ c= - \theta$, $ \Delta_c=  ]-\infty, 0[$ and we find a new regime in the asymptotic expansions of $a_T$, $A_T$, $B_T$.

\begin{lem} \label{L-AT4}
For $ c= - \theta$, there exists a unique $a_T$, which belongs to the interior of $\Delta_c$ and converges to the origin 
as $T$ goes to infinity, solution of the implicit equation
\begin{equation}
\label{IMPEQNN}
\cL^{\prime}(a) + \frac{1}{T}\cH^{\prime}(a) = 0
\end{equation}
where the functions $\cL$ and $\cH$ are given by \eqref{E-DEFL} and \eqref{E-DEFH}.
Moreover, there exists a sequence $(\gamma_k)$ such that, for any $p > 0$ and $T$ large enough,
\begin{equation}
A_T = \exp \left(- T I(c) \right) (e \theta T)^{1/4} \left[ 1 + \sum_{k=1}^{2p} \frac{\gamma_k}{(\sqrt{T})^k} + {\mathcal O}\Bigl( \frac{1}{T^{p} \sqrt{T}}\Bigr) \right]. \label{DEVAT4FIN}
\end{equation}
The sequence $(\gamma_k)$ only depends on the Taylor expansion of $a_T$ at the neighborhood of the origin together with the derivatives of $\cL$ and
$\cH$ at $0$. For example,
\begin{equation*}
\gamma_1 = \frac{3}{8 \sqrt{\theta}}.
\end{equation*}
\end{lem}

\noindent {\bf Proof.} The proof of Lemma \ref{L-AT4} is given in Appendix B.\demend

\noindent
It now remains to give the expansion for $B_T$.

\begin{lem}$\pt$ \label{L-BT4}
For $ c= - \theta$, there exists a sequence $(\delta_k)$ such that, for any $p>0$ and $T$ large enough,
\begin{equation} \label{DEVBT4}
B_T= \sum_{k=1}^{2p} \frac{\delta_k}{(\sqrt{T})^k} + {\mathcal O}\Bigl( \frac{1}{T^{p} \sqrt{T}}\Bigr).
\end{equation}
The sequence $(\delta_k)$ only depends on the Taylor expansion of $a_T$ at the neighborhood of the origin together with the derivatives of $\cL$ and
$\cH$ at $0$. For example,
\begin{equation*}
\delta_1=\frac{1}{2 \pi} e^{-1/4} \Gamma\left( \frac{1}{4} \right).
\end{equation*}
\end{lem}

\noindent
{\bf Proof.} The proof of Lemma \ref{L-BT4} is given in Appendix C.\demend

\noindent
{\bf Proof of Theorem \ref{T-SLDPE} d).} We immediately deduce \eqref{SLDPE4}  from \eqref{DEVAT4FIN} together with \eqref{DEVBT4}.
\demend

\subsection{Proof of Theorem \ref{T-SLDPE} e).}

We obtain from \eqref{SOLOUP} that $X_T$ is Gaussian with $\cN(0, \sigma_T^2)$ distribution where
\begin{equation*}
\sigma_T^2=\frac{1}{2\theta}\Bigl(\exp(2 \theta T) - 1 \Bigr).
\end{equation*}
Moreover, we clearly have
\begin{eqnarray}
\p(\widehat{\theta}_T \leq 0)&=&\p(X_T^2 \leq T)=\p(|X_T| \leq \sqrt{T})=2\p(0 \leq X_T \leq \sqrt{T}) \nonumber \\
&=& 2\p(0 \leq Z \leq d_T)
\label{DEVITO}
\end{eqnarray}
where $Z$ is an $\cN(0, 1)$ random variable and the sequence  $(d_T)$ satisfies
\begin{equation*}
d_T= \sqrt{2 \theta T \exp(-2 \theta T)}\Bigl[1 + {\mathcal O} (\exp(-2 \theta T)) \Bigr].
\end{equation*}
For all $x>0$, denote
$$
F(x)=\int_0^x f(t)\,dt 
$$
where $f$ stands for the probability density function of the $\cN(0, 1)$ distribution. It is well-known that
for all $n \geq 1$, the Gaussian derivatives
$$
f^{(n)}(x)= \frac{(-1)^n}{2^{n/2}}H_n\Bigl(\frac{x}{\sqrt{2}}\Bigr) f(x)
$$
where $(H_n)$ is the sequence of Hermite polynomials. For example, we have $H_0(x)=1$, $H_1(x)=2x$, $H_2(x)=-2+4x^2$, etc.
Hence,
$$
f^{(n)}(0)= \frac{(-1)^n}{\sqrt{2 \pi}2^{n/2}}H_n
$$
where $H_n=H_n(0)$ are the Hermite numbers given by the recurrence relation $H_n=-2(n-1)H_{n-2}$ with
$H_0=1$ and $H_1=0$ which implies that
\begin{eqnarray*}
H_n= \left \{ \begin{array}{lll}
    0 & \text{ if }  n \text{ is odd, } \vspace{1ex}\\
    {\displaystyle \hspace{-1ex} \frac{(-1)^{n/2} n!}{(n/2)!} } & \text{ if }  n \text{ is even. } 
\end{array}  \right.
\end{eqnarray*}
Consequently, at the neighborhood of the origin, we have for all $x>0$ the Taylor expansion 
\begin{equation}
\label{TAYLORF}
F(x)=\frac{x}{\sqrt{2 \pi }}
\left[1 + \sum_{k=1}^{p}\frac{(-1)^{k} x^{2k}}{(2k+1)2^k k!}
+ {\cal O}(x^{p+1})\right].
\end{equation}
Therefore, we deduce from the identity $\p(\widehat{\theta}_T \leq 0)=2F(d_T)$ together with
\eqref{TAYLORF} that
\begin{equation*}
\p(\widehat{\theta}_T \leq 0) = 2
\frac{\exp(-\theta T) \sqrt{2 \theta T}}{\sqrt{2 \pi}}
\left[1 + \sum_{k=1}^{p}\frac{(-1)^{k} (\theta T e^{-2\theta T})^k}{(2k+1) k!} 
+ {\cal O}\Bigl((Te^{-2\theta T})^{p+1}\Bigr)\right],
\end{equation*}
which immediately leads to \eqref{SLDPE5}.

\demend

\renewcommand{\thesection}{\Alph{section}}
\renewcommand{\theequation}{\thesection.\arabic{equation}}
\setcounter{section}{0}
\setcounter{subsection}{0}
\setcounter{equation}{0}
\section{Appendix A: Proof of the keystone Lemma \ref{L-LEMFOND}.}

Our goal is to prove the asymptotic expansion \eqref{E-DECFOND} associated to the normalized cumulant generating function
$\cL_T$. Via the same approach as in Section 17.3 of Liptser and Shiryaev \cite{LS01}, we have
\begin{eqnarray*}
\cL_T(a) &=& \frac{1}{T} \log \E \! \left[ \exp\Bigl(a \int_0^T X_t dX_t - ac \int_0^T X_t^2 \,dt\Bigr) \right], \\
&=& \frac{1}{T} \log \E_{\vp} \! \!\left[ \exp\Bigl((a+\theta - \vp) \int_0^T X_t dX_t +\frac{1}{2}( -2ac  - \theta^2 + \vp^2) \int_0^T X_t^2 \,dt\Bigr) \right]
\end{eqnarray*}
for all $\vp \in \R$, where $\E_{\vp}$ stands for the expectation after the change of measures
$$
\frac{d \p_{\vp}}{d \p}= 
\exp \Bigl( ( \vp- \theta) \int_0^T X_t dX_t - \frac{1}{2} (\vp^2 - \theta^2) \int_0^T X_t^2\,dt\Bigr).
$$
Hereafter, consider $a \in \Delta_c=\{ a \in \R,  ~~\theta^2 +2ac >0,~~ a + \theta < \sqrt{ \theta^2+2ac} \}$ so that we can choose 
$\vp=\vp(a)$ where
$
\vp(a)= -\sqrt{\theta^2 +2ac}.
$
Then, if we denote $ \tau(a) = a + \theta - \vp(a)$, we obtain that
\begin{equation}
\label{DECOMGF1}
\cL_T(a)= \frac{1}{T} \log \E_{\vp} \! \!\left[ \exp\Bigl(\tau(a) \int_0^T X_t dX_t \Bigr) \right].
\end{equation}
However, we have from It\^o's formula that
\begin{equation*}
\int_0^{T} X_t dX_t= \frac{1}{2} \Bigl( X_{T}^2-T\Bigr).
\end{equation*}
Consequently, we obtain from \eqref{DECOMGF1} that
\begin{equation}
\label{DECOMGF2}
\cL_T(a)= -\frac{\tau(a)}{2}+ \frac{1}{T} \log \E_{\vp} \! \!\left[ \exp\Bigl(\frac{\tau(a)}{2} X_T^2 \Bigr) \right].
\end{equation}
Under the measure $\p_{\vp}$, $X_T$ is a Gaussian random variable with zero mean and variance $\sigma^2_T(a)$ given by 
$$
\sigma^2_T(a) = -\frac{1 - \exp(2 \vp(a) T) }{2 \vp(a)}.
$$
This variance, together with \eqref{DECOMGF2}, leads to
\begin{equation}
\label{DECOMGF3}
\cL_T(a)= -\frac{\tau(a)}{2}- \frac{1}{2T} \log \left( 1 + \frac{\tau(a)}{2 \vp(a)} \Bigl(1-\exp(2 \vp(a) T)  \Bigr) \right).
\end{equation}
Finally, if
$$
h(a)= \frac{a + \theta}{\vp(a)} = \frac{\tau(a)}{\vp(a)} +1,
$$
we find from \eqref{DECOMGF3} the decomposition
\begin{eqnarray*}
\cL_T(a)
&=& -\frac{\tau(a) }{2} - \frac{1}{2T} \log \left( 1 + \frac{1}{2}(h(a) -1)\Bigl(1 - \exp(2\vp(a) T)\Bigr) \right),\\
&=& -\frac{\tau(a)}{2}  - \frac{1}{2T} \log \left( \frac{1}{2}(1 + h(a)) +  \frac{1}{2}(1 - h(a)) \exp(2\vp(a) T) \right),\\
&=& -\frac{\tau(a)}{2}  - \frac{1}{2T} \log \left( \frac{1}{2}(1 + h(a)) \right) 
- \frac{1}{2T} \log \left( 1+   \frac{1 - h(a)}{1 + h(a)} \exp(2\vp(a) T) \right),\\
&=& \cL(a) + \frac{1}{T} \cH(a) + \frac{1}{T} \cR_T(a).
\end{eqnarray*}
One can observe that the remainder $\cR_T(a)$ goes to zero exponentially fast as
$\cR_T(a) = \cO ( \exp(2 \varphi(a) T) )$, which completes the proof of Lemma \ref{L-LEMFOND}. \demend

\section{Appendix B: On the expansions of $A_T$.}
\setcounter{equation}{0}

All asymptotic expansions associated with $A_T$ are related on the fact that there exists a unique $a_T$, 
which belongs to the interior of $\Delta_c$ and converges to its border 
$a_c=2(c - \theta)$ if $c> \theta$, and to the origin if $[c|< \theta$, solution of the implicit equation
\begin{equation}
\label{GENIMPLEQ}
\cL^{\prime}(a) + \frac{1}{T}\cH^{\prime}(a) = 0
\end{equation}
where the functions $\cL$ and $\cH$ are given by \eqref{E-DEFL} and \eqref{E-DEFH}.
After some straightforward calculation, \eqref{GENIMPLEQ} can be rewritten as
\begin{equation} 
\label{EQDEFAT}
T \vp(a) (\vp(a) - c) (\vp(a) + a + \theta) = c(a + \theta) - \vp^2(a).
\end{equation}
One can observe that \eqref{EQDEFAT} may be rewritten as
\begin{equation*} 
T \vp(a) (\vp(a) - c) (\vp(a) + \theta) ( \vp(a) +2c - \theta)= -\frac{c}{2}(\vp^2(a)+\theta^2 - 2 \theta c)
\end{equation*}
which ensures that $\varphi(a_T)$ converges to $\theta - 2c$, while $a_T<2(c - \theta)$
and converges to $a_c$.
Moreover, it follows from \eqref{E-DECFOND} that 
\begin{eqnarray} 
A_T &=& \exp \Bigl( T \cL(a_T) + \cH(a_T) + \cR_T(a_T) \Bigr),\notag\\
	&=& \exp \Bigl( T \cL(a_T) \Bigr) \exp \Bigl( \cH(a_T) \Bigr) \exp \Bigl( \cR_T(a_T) \Bigr).
\label{DEVAT2}
\end{eqnarray}
Therefore, the proofs of the expansions of $A_T$ are divided into four steps. 
The first one is devoted to the asymptotic expansions of $a_T$ and $\vp(a_T)$. The last three one
deal with the asymptotic expansions of all terms in \eqref{DEVAT2}. 

\subsection{Proof of Lemma \ref{L-AT2}.}

{\bf Step 1.} {\it One can find two sequences $(a_k)$ and $(\vp_k)$ such that, for any $p>0$ and $T$ large enough,
\begin{equation*} 
a_T =\sum_{k=0}^p \frac{a_k}{T^k}+{\cal O}\Bigl(\frac{1}{T^{p+1}}\Bigr)
\end{equation*} 
where $a_0 = 2(c - \theta)$,
\begin{equation*} 
a_1 = \frac{\theta -2c}{ 3c - \theta} \hspace{1cm} \text{and} \hspace{1cm} 
a_2 =  -\frac{c(c^2 -5 \theta c +2\theta^2)}{2(c - \theta)(3c-\theta)^3},
\end{equation*} 
 \begin{equation*} 
\varphi(a_T) = \sum_{k=0}^p  \frac{\varphi_k}{T^k}+{\cal O}\Bigl(\frac{1}{T^{p+1}}\Bigr)
\end{equation*} 
where  $\vp_0 = \theta - 2c$,
\begin{equation*} 
\vp_1 = \frac{c}{ 3c - \theta}
\hspace{1cm} \text{and} \hspace{1cm} 
\varphi_2 = \frac{c^2(4c^2 - 9 \theta c +3 \theta^2)}{2(c-\theta)(2c-\theta)(3c -\theta)^3}.
\end{equation*} 
}
\noindent
{\bf Proof.} We are in the situation where
$c > \theta$, $a_c = 2(c - \theta)$ and $\vp(a_c) = \theta - 2c$. Consequently, $\vp(a_c) - c =  \theta - 3c \ne 0$ while $\vp(a_c) + a_c + \theta = 0$. One can easily deduce from \eqref{EQDEFAT} that
\begin{equation} \label{JUSTIF}
\lim_{T \to \infty} T (\vp(a_T) + a_T + \theta) = \frac{c-\theta}{\theta - 3c}.
\end{equation}
Therefore, the conjunction of \eqref{EQDEFAT} and \eqref{JUSTIF} leads to
the asymptotic expansions of $a_T$ and $\vp(a_T)$. Let us show how to calculate the first terms of the expansions. On the one hand, as 
$$\vp(a_T)= - \sqrt{\theta^2 +2 a_T c},$$
we have
\begin{equation} \label{LINKAPHI}
a_0 = \frac{\vp_0^2 - \theta^2}{2c}, \qquad a_1 = \frac{\vp_0 \vp_1}{c}, \qquad a_2 = \frac{2\vp_0 \vp_2 + \vp_1^2}{2c}.
\end{equation}
On the other hand, it follows from \eqref{EQDEFAT} that
\begin{eqnarray*}
\vp_1 + a_1 &=& \frac{c(a_0+\theta) - \vp_0^2}{\vp_0(\vp_0 -c)},\\
\vp_2 + a_2 &= &\frac{c a_1 - 2\vp_0 \vp_1 -(\vp_1 + a_1)\vp_1(2\vp_0  - c)}{\vp_0(\vp_0-c)}.
\end{eqnarray*}
Finally, in order to calculate $a_1, a_2, \vp_1$, and $\vp_2$, it is only necessary to solve very simple linear
systems. The rest of the proof is left to the reader.
\demend
\ \vspace{-1ex} \\
{\bf Step 2.} {\it  One can find a sequence $(\alpha_k)$ such that, for any $p > 0$ and $T$ large enough,
\begin{equation} 
\label{DEVLT2}
\exp \left( T \cL(a_T) \right) = \exp \left( -T I(c) + \frac{1}{2} \right) \left[ 1 + \sum_{k=1}^p \frac{\alpha_k}{T^k} + {\mathcal O}\Bigl( \frac{1}{T^{p+1}}\Bigr) \right].
\end{equation}
The sequence $(\alpha_k)$ only depends on $(a_k)$ together with the derivatives of $\cL$ at point $a_c$. For example,
\begin{equation*} 
\alpha_1 = \frac{c(c^2-3\theta c +\theta^2)}{2(c - \theta)(2c-\theta)(3c-\theta)^2}.
\end{equation*}
}
\noindent
{\bf Proof.} By the Taylor expansion of $\cL$ at the neighborhood of $a_c$, we have the existence 
of a sequence $(\ell_k)$ such that, for any $p > 0$ and $T$ large enough,
\begin{equation}
\label{DEVLTT}
T \cL(a_T) = T \cL(a_c) + a_1 \cL^{\prime}(a_c) + \sum_{k=1}^p \frac{\ell_k}{T^k} + {\mathcal O}\Bigl( \frac{1}{T^{p+1}}\Bigr).
\end{equation}
On the one hand,
$$
a_1 \cL^{\prime}(a_c) = \frac{1}{2}.
$$
On the other hand,
$$
\ell_1 = a_2 \cL^{\prime}(a_c) + \frac{1}{2} a_1^2 \cL^{\prime \prime}(a_c)= \frac{c(c^2-3\theta c +\theta^2)}{2(c - \theta)(2c- \theta)(3c-\theta)^2}.
$$
Therefore, \eqref{DEVLT2} clearly follows from \eqref{DEVLTT}.
\demend
\ \vspace{-1ex} \\
{\bf Step 3.} {\it One can find a sequence $(\beta_k)$ such that, for any $p > 0$ and $T$ large enough,
\begin{equation} \label{DEVHT2}
\exp \left( \cH(a_T) \right) = \sqrt{ \frac{ 2 \vp_0 T }{ \vp_1 + a_1 } } \left[ 1 + \sum_{k=1}^p \frac{\beta_k}{T^k} + {\mathcal O}\Bigl( \frac{1}{T^{p+1}}\Bigr) \right].
\end{equation}
The sequence $(\beta_k)$ only depends on $(a_k)$ together with the derivatives of $\cH$ at point $a_c$. For example,
\begin{equation*}
\beta_1 = \frac{c(c^2-3\theta c +\theta^2)}{(c - \theta)( \theta - 2c)(3c-\theta)^2}.
\end{equation*}
}
\noindent
{\bf Proof.} By the very definition of $\cH$, we have
$$
\exp \left( \cH(a_T) \right) = \sqrt{ \frac{ 2 \vp(a_T) T }{ T ( \vp(a_T) + a_T + \theta ) }}.
$$
Consequently, the expansion of the square root together with those of $a_T$ and $\vp(a_T)$ ensure the existence of a sequence $(\beta_k)$ such that \eqref{DEVHT2} is true. Moreover, as for $(\alpha_k)$, the sequence $(\beta_k)$ can be explicitly calculated. For example
$$
\beta_1 =  \frac{1}{2} \left( \frac{\vp_1}{\vp_0} - \frac{\vp_2 + a_2}{\vp_1 + a_1} \right)= \frac{c(c^2-3\theta c +\theta^2)}{(c - \theta)( \theta - 2c)(3c-\theta)^2}.
$$
\vspace{-1ex} \\
{\bf Step 4.} {\it The remainder $\cR_T(a_T)$ goes to zero exponentially fast
\begin{equation}
\label{DEVRT2}
\cR_T(a_T) = \cO \Bigl(T \exp(2 \vp_0 T) \Bigr).
\end{equation}
}
\noindent
{\bf Proof.} 
The result follows from \eqref{E-DEFR}  together with the fact that
$\vp_0<-\theta<0$. More precisely, we have
$$
\frac{1 - h(a_T)}{ 1 + h(a_T)} 
= \frac{T (\vp(a_T) - a_T - \theta)}{T (\vp(a_T) + a_T + \theta)}
$$
which implies via \eqref{JUSTIF} that 
\begin{equation}
\label{LIMRTT}
\lim_{T \to \infty} \frac{1}{T}\Bigl(\frac{1 - h(a_T)}{ 1 + h(a_T)}\Bigr) = \frac{\vp_0 - a_0 - \theta}{\vp_1 + a_1} = \frac{2(2c-\theta)(3c-\theta)}{c-\theta}.
\end{equation}
Consequently, we immediately deduce \eqref{DEVRT2} from \eqref{E-DEFR} and \eqref{LIMRTT}. \demend
\ \vspace{-1ex} \\
{\bf Proof of Lemma \ref{L-AT2}.}
It follows from the conjunction of \eqref{DEVAT2}, \eqref{DEVLT2}, \eqref{DEVHT2} and \eqref{DEVRT2}
that there exists a sequence $(\gamma_k)$ such that, for any $p > 0$ and $T$ large enough,
\begin{eqnarray}
A_T &=& \exp \left(- T I(c) + \frac{1}{2} \right) \sqrt{ \frac{ 2 \vp_0 T }{ \vp_1 + a_1 } } \left[ 1 + \sum_{k=1}^p \frac{\gamma_k}{T^k} + {\mathcal O}\Bigl( \frac{1}{T^{p+1}}\Bigr) \right],\notag\\
&=& \exp \left(- T I(c) + P(c) \right) \sqrt{e T } \left[ 1 + \sum_{k=1}^p \frac{\gamma_k}{T^k} + {\mathcal O}\Bigl( \frac{1}{T^{p+1}}\Bigr) \right], \label{DEVATFIN}
\end{eqnarray}
where $P(c)$ is given by \eqref{DEFPQ1}. Finally, the sequence $(\gamma_k)$ can be explicitly calculated  by use of $(a_k)$
together with the derivatives of $\cL$ and $\cH$ at point $a_c$. For example,
$$
\gamma_1 =  \alpha_1 + \beta_1= \frac{c(c^2-3\theta c +\theta^2)}{2(c - \theta)( \theta - 2c)(3c-\theta)^2}.
$$
\demend

\subsection{Proof of Lemma \ref{L-AT3}.}

We are in the situation where $|c| < \theta$ and $c \ne 0$ which means that $a_c = 0$ and $\vp(a_c) = -\theta$. Consequently, 
$\vp(a_c) - c = -(\theta + c) \ne 0$ while $\vp(a_c) + a_c + \theta = 0$.  The proof of Lemma \ref{L-AT3}
follows exactly the same lines as those of Lemma \ref{L-AT2}. The only notable thing to mention is that
\begin{equation*} 
a_T =\sum_{k=0}^p \frac{a_k}{T^k}+{\cal O}\Bigl(\frac{1}{T^{p+1}}\Bigr)
\end{equation*} 
where $a_0 = 0$,
\begin{equation*} 
a_1 =  -\frac{\theta}{c + \theta}\hspace{1cm} \text{and} \hspace{1cm} 
a_2 =  -\frac{c(c^2 + 3 \theta c - 2 \theta^2)}{2(c - \theta)(c + \theta)^3}, 
\end{equation*} 
\begin{equation*} 
\varphi(a_T) = \sum_{k=0}^p  \frac{\varphi_k}{T^k}+{\cal O}\Bigl(\frac{1}{T^{p+1}}\Bigr)
\end{equation*} 
where $\vp_0 = -\theta$,
\begin{equation*} 
\vp_1 = \frac{c}{c + \theta}\hspace{1cm} \text{and} \hspace{1cm} 
\varphi_2 = \frac{c^2(2c^2 + 3\theta c - 3\theta^2)}{2\theta (c - \theta)(c + \theta)^3}. 
\end{equation*}
Therefore, the rest of the proof of the Lemma \ref{L-AT3} is left to the reader. \demend

\subsection{Proof of Lemma \ref{L-AT4}.}

The proof of Lemma \ref{L-AT4} is slightly different from the one of Lemma \ref{L-AT2}. More precisely, there is
a change of regime in the asymptotic expansions of $a_T$ and $\vp(a_T)$. 
\ \vspace{1ex} \\
{\bf Step 1.} {\it One can find two sequences $(a_k)$ and $(\vp_k)$ such that, for any $p>0$ and $T$ large enough,
\begin{equation*} 
a_T = \sum_{k=0}^{2p} \frac{a_k}{(\sqrt{T})^k}+{\mathcal O}\Bigl( \frac{1}{T^{p} \sqrt{T}}\Bigr)
\end{equation*} 
where $a_0 = 0$, $a_1 = - \sqrt{\theta}$, and $a_2=-1/8$,
\begin{equation*} 
\varphi(a_T) = \sum_{k=0}^{2p} \frac{\vp_k}{(\sqrt{T})^k}+{\mathcal O}\Bigl( \frac{1}{T^{p} \sqrt{T}}\Bigr)
\end{equation*} 
where  $\vp_0 = -\theta$, $\vp_1 = - \sqrt{\theta}$, and $\vp_2=3/8$.
}
\ \vspace{1ex} \\
{\bf Proof.} We are in the situation where $c=-\theta$, $a_c = 0$ and $\vp(a_c) = -\theta$ which clearly implies that
$\vp(a_c) - c = -(\theta + c) = 0$ and $\vp(a_c) + a_c + \theta = 0$. It leads to a change of regime 
in the expansions of $a_T$ and $\vp(a_T)$ comparing to the expansions of $a_T$ and $\vp(a_T)$
in Lemma \ref{L-AT2}. As a matter of fact, we obtain from \eqref{EQDEFAT} that
\begin{equation} \label{JUSTIF2}
\lim_{T \to \infty} T (\vp(a_T) + \theta)(\vp(a_T) + a_T + \theta) = 2\theta.
\end{equation}
Therefore, one can easily deduce the expansions of $a_T$ and $\vp(a_T)$
from \eqref{EQDEFAT} and \eqref{JUSTIF2}. The calculation of the first terms is straightforward. For example, as 
$$\vp(a_T)= - \sqrt{\theta^2 -2\theta a_T },$$ 
we obtain that $a_0=0$, 
$\vp_0=-\theta$,
\begin{equation*}
\vp_1=a_1 \hspace{1cm} \text{and} \hspace{1cm} \vp_2=a_2+\frac{1}{2}.
\end{equation*}
In addition, we infer from \eqref{EQDEFAT} that
\begin{equation*}
\vp_1(a_1+\vp_1)= 2 \theta \hspace{1cm} \text{and} \hspace{1cm} 
a_2 + 3 \vp_2=1.
\end{equation*}
Consequently, we immediately obtain that $a_1^2=\theta$ which implies that $a_1=-\sqrt{\theta}$ 
as $a_T$ belongs to the interior of $\Delta_c=  ]-\infty, 0[$. It remains to solve the simple linear
system
\begin{equation*}
\begin{cases}
2a_2 - 2 \vp_2 =-1,\vspace{1ex}\\
a_2+3 \vp_2=1
\end{cases}
\end{equation*}
which solution is $a_2=-1/8$ and $\vp_2=3/8$. \demend
\ \vspace{-1ex} \\
{\bf Step 2.} {\it  One can find a sequence $(\alpha_k)$ such that, for any $p > 0$ and $T$ large enough,
\begin{equation} 
\label{DEVLT3}
\exp \left( T \cL(a_T) \right) = \exp \left( -T I(c) + \frac{1}{4} \right) \left[ 1 + \sum_{k=1}^{2p} \frac{\alpha_k}{(\sqrt{T})^k} + {\mathcal O}\Bigl( \frac{1}{T^{p} \sqrt{T}}\Bigr) \right].
\end{equation}
The sequence $(\alpha_k)$ only depends on $(a_k)$ together with the derivatives of $\cL$ at the origin. For example,
\begin{equation*} 
\alpha_1 = -\frac{3}{16 \sqrt{\theta}}.
\end{equation*}
}
\noindent
{\bf Proof.} By the Taylor expansion of $\cL$ at the neighborhood of the origin, as $\cL^{\prime}(0)=0$,
we have the existence of a sequence $(\ell_k)$ such that, for any $p > 0$ and $T$ large enough,
\begin{equation}
\label{DEVLTTT}
T \cL(a_T) = T\cL(0) + \frac{a_1^2}{2} \cL^{(2)}(0)  
 +  \sum_{k=1}^{2p} \frac{\ell_k}{(\sqrt{T})^k} + {\mathcal O}\Bigl( \frac{1}{T^{p} \sqrt{T}}\Bigr).
\end{equation}
On the one hand, $\cL^{(2)}(0)=1/(2\theta)$ which implies that
$$
\frac{a_1^2}{2} \cL^{(2)}(0) = \frac{1}{4}.
$$
On the other hand, as $\cL^{(3)}(0)=3/(2\theta^2)$, we also have
$$
\ell_1 = a_1 a_2 \cL^{(2)}(0) + \frac{a_1^3}{6} \cL^{(3)}(0)=-\frac{3}{16 \sqrt{\theta}}.
$$
Therefore, we deduce \eqref{DEVLT3} from \eqref{DEVLTTT}.
\demend
\ \vspace{-1ex} \\
{\bf Step 3.} {\it One can find a sequence $(\beta_k)$ such that, for any $p > 0$ and $T$ large enough,
\begin{equation} \label{DEVHT3}
\exp \left( \cH(a_T) \right) = (\theta T )^{1/4} \left[ 1 + \sum_{k=1}^{2p} \frac{\beta_k}{(\sqrt{T})^k} + {\mathcal O}\Bigl( \frac{1}{T^{p} \sqrt{T}}\Bigr) \right].
\end{equation}
The sequence $(\beta_k)$ only depends on $(a_k)$ together with the derivatives of $\cH$ at the origin. For example,
\begin{equation*}
\beta_1 = \frac{9}{16 \sqrt{\theta}}.
\end{equation*}
}
\noindent
{\bf Proof.} By the very definition of $\cH$, we have
$$
\exp \left( \cH(a_T) \right) = \sqrt{ \frac{ 2 \vp(a_T) \sqrt{T} }{ \sqrt{T} ( \vp(a_T) + a_T + \theta ) }}.
$$
Hence, the expansion of the square root together with those of $a_T$ and $\vp(a_T)$ ensure the existence of a sequence $(\beta_k)$ such that \eqref{DEVHT3} is true. As before, the sequence $(\beta_k)$ can be explicitly calculated. For example, 
$$
\beta_1 =  \frac{1}{2} \left( \frac{\vp_1}{\vp_0} - \frac{\vp_2 + a_2}{\vp_1 + a_1} \right) = \frac{9}{16 \sqrt{\theta}}.
$$
\vspace{-1ex} \\
{\bf Step 4.} {\it The remainder $\cR_T(a_T)$ goes to zero exponentially fast
\begin{equation}
\label{DEVRT3}
\cR_T(a_T) = \cO \Bigl(\sqrt{T} \exp(-2 \theta T) \Bigr).
\end{equation}
}
\noindent
{\bf Proof.}  We have
$$
\frac{1 - h(a_T)}{ 1 + h(a_T)} 
= \frac{\sqrt{T} (\vp(a_T) - a_T - \theta)}{\sqrt{T} (\vp(a_T) + a_T + \theta)},
$$
which implies that 
\begin{equation}
\label{LIMRTTT}
\lim_{T \to \infty} \frac{1}{\sqrt{T}}\Bigl(\frac{1 - h(a_T)}{ 1 + h(a_T)}\Bigr) = \sqrt{\theta}.
\end{equation}
Consequently, we immediately deduce \eqref{DEVRT3} from \eqref{E-DEFR} and \eqref{LIMRTTT}. \demend
\ \vspace{-1ex} \\
{\bf Proof of Lemma \ref{L-AT4}.}
It follows from \eqref{DEVAT2} together with \eqref{DEVLT3}, \eqref{DEVHT3} and \eqref{DEVRT3}
that there exists a sequence $(\gamma_k)$ such that, for any $p > 0$ and $T$ large enough,
\begin{equation*}
 \exp \left(- T I(c) \right) (e \theta T)^{1/4} \left[ 1 + \sum_{k=1}^{2p} \frac{\gamma_k}{(\sqrt{T})^k} + {\mathcal O}\Bigl( \frac{1}{T^{p} \sqrt{T}}\Bigr) \right]
\end{equation*}
where the sequence $(\gamma_k)$ can be explicitly calculated  by use of $(a_k)$
together with the derivatives of $\cL$ and $\cH$ at the origin. For example,
$$
\gamma_1 =  \alpha_1 + \beta_1= \frac{3}{8\sqrt{\theta}}.
$$
\demend

\section{Appendix C: On the expansions of $B_T$.}
\setcounter{equation}{0}

\subsection{General considerations.}
In order to unify the notations, let $\alpha_T= a_c$ if $c < -\theta$ and $\alpha_T= a_T$ otherwise. In addition, denote
\begin{eqnarray*}
\beta_T= \left \{ \begin{array}{lll}
    \sigma_c \sqrt{T} & \text{ if } \hspace{2ex}  c < -\theta, \vspace{1ex}\\
    \sqrt{T}  & \text{ if } \hspace{2ex} c =-\theta, \vspace{1ex}\\
     \hspace{1ex} T   & \text{ if } \hspace{2ex} |c| < \theta, \vspace{1ex}\\
      - T  &\text{ if } \hspace{2ex} c > \theta.
\end{array}  \right.
\end{eqnarray*}
One can observe that we always have $\alpha_T \beta_T < 0$. Then, in all different cases, 
\begin{equation} 
\label{BT}
B_T = \E_T\Bigl[\exp(-\alpha_T \beta_T U_T) \ind_{U_T \leq 0}\Bigr]
\end{equation}
where
\begin{equation*}
U_T = \frac{Z_T(c)}{\beta_T}.
\end{equation*}
Denote by $\Phi_T$ the characteristic function of $U_T$ under $\p_T$ and assume in all the sequel that $c \neq 0$.

\begin{lem} \label{L-GENE}
For $T$ large enough, $\Phi_T$ belongs to $\dL^2(\R)$ and, for all $u \in \R$,
\begin{equation} 
\label{PHITOLT}
\Phi_T(u) = \exp \left( T \cL_T \left( \alpha_T +  \frac{i u}{\beta_T} \right)- T \cL_T(\alpha_T) \right).
\end{equation}
Moreover, we can split $B_T$ into two terms, $B_T=C_T+D_T$ where
\begin{eqnarray}
C_T &=& -\frac{1}{2\pi \alpha_T \beta_T} \int_{|u| \leq s_T}\left(1+\frac{iu}{\alpha_T \beta_T}\right)^{-1}\Phi_{T}(u) du, \label{CT}\\
D_T &=& -\frac{1}{2\pi \alpha_T \beta_T}  \int_{|u| > s_T}\left(1+\frac{iu}{\alpha_T \beta_T}\right)^{-1}\Phi_{T}(u) du. \label{DT}
\end{eqnarray}
where $s_T$ is chosen in such a way that there are positive constants $C$ and $0<\nu<1$ satisfying
\begin{equation} 
\label{HYPHYP}
\min \left( \frac{T s_T^2}{\beta_T^2}, \frac{T \sqrt{s_T}}{\sqrt{|\beta_T|}} \right) \geq C T^{\nu}
\end{equation}
and there exist two positive constants $d$ and $D$ such that
\begin{equation} \label{MAJDT}
|D_T|\leq d~T~\exp(-D T^{\nu}).
\end{equation}
\end{lem}

We choose $s_T$ large enough to satisfy \eqref{HYPHYP} and small enough to enable us to intervene integral and summation into \eqref{CT}. 
The expansion of $C_T$ thus follows from that of $\Phi_T$ and some tedious calculations. Finally, 
\eqref{MAJDT} tells us that the expansion of $B_T$ is nothing but that of $C_T$.\\

\noindent{\bf Proof of Lemma \ref{L-GENE}.}
For all $u \in \R$, we have
\begin{eqnarray*}
\Phi_T(u) &=& \E_T\Bigl[\exp(iu U_T)\Bigr],\\
&=&\es{\exp \left( i u \frac{Z_T(c)}{\beta_T} \right) \exp \left( \alpha_T Z_T(c) - T \cL_T \alpha_T ) \right)},\\
&=&\es{\exp \left(\left( \alpha_T  +  \frac{i u}{\beta_T} \right) Z_T(c)\right)} \exp \left( - T \cL_T(\alpha_T) \right),\\
&=&\exp \left( T \cL_T \left( \alpha_T +  \frac{i u}{\beta_T} \right)- T \cL_T(\alpha_T) \right).
\end{eqnarray*}
We shall see in Appendix D that for $T$ large enough, $\Phi_T \in \dL^2(\R)$. Then, it follows from Parseval formula that
\begin{eqnarray*}
B_T &= & \E_T\Bigl[ \exp(-\alpha_T \beta_T U_T) \ind_{U_T \leq 0}\Bigr],\\
&=&-\frac{1}{2 \pi} \int_{\R} \frac{1}{\alpha_T \beta_T +i u }~ \Phi_T(u) du, \\
&=&-\frac{1}{2 \pi \alpha_T \beta_T} \int_{\R} \left( 1 + \frac{iu}{\alpha_T \beta_T} \right)^{-1} \Phi_{T}(u) du
\end{eqnarray*}
which implies that $B_T=C_T+D_T$ where $C_T$ and $D_T$ are given by \eqref{CT} and \eqref{DT}.
It remains to show that $D_T$ goes exponentially fast to zero. We deduce from Cauchy-Schwarz inequality that
\begin{equation}
\label{MAJDTCS}
|D_T|^2 \leq \frac{1}{4 \pi^2 \alpha_T^2 \beta_T^2}
 \int_{|u| > s_T} \left( 1 + \frac{u^2}{ (\alpha_T \beta_T)^2 } \right)^{-1} \, du
\int_{|u| > s_T} |\Phi_T(u) |^2 \, du.
\end{equation}
On the one hand, 
\begin{equation} \label{INTDT}
\int_{|u| > s_T} \left( 1 + \frac{u^2}{ (\alpha_T \beta_T)^2 } \right)^{-1} \, du \leq 
|\alpha_T \beta_T| \int_{\R} \frac{1}{1 + v^2} \, dv \leq
|\alpha_T \beta_T| \pi.
\end{equation}
On the other hand, we deduce from \eqref{PHITOLT} together with inequality \eqref{MAJPHIT} that for $T$ large enough,
$$
|\Phi_T(u)|^2 \leq 4
\ell(\alpha_T,c,\theta) \Bigl( 1 + \gamma_T^2 u^2 \Bigr)^{1/4}
\exp \left( \frac{T\varphi_T}{8} \gamma_T^2  u^2 \Bigl( 1 + \gamma_T^2 u^2
\Bigr)^{-3/4} \right)
$$
where $\varphi_T=\varphi(\alpha_T) $ and
$$
\gamma_T= \frac{2|c|}{|\beta_T| \varphi^2(\alpha_T)}.
$$
It is not hard to see that we can find a positive constant $C_{\ell}$ such that, for $T$ large enough,
$\ell(\alpha_T,c,\theta) \leq C_{\ell}  T. $
Consequently, if $\delta_T=\gamma_T s_T$, we obtain that
\begin{eqnarray} 
\int_{|u| > s_T} |\Phi_T(u) |^2 \, du &\leq &
8 C_{\ell}  T \int_{s_T}^\infty 
\Bigl( 1 + \gamma_T^2 u^2 \Bigr)^{1/4}
\exp \left( \frac{T\varphi_T}{8} \gamma_T^2  u^2 \Bigl( 1 + \gamma_T^2 u^2
\Bigr)^{-3/4} \right)\,du,\notag\\
&\leq &
\frac{ 8 C_{\ell}  T}{\gamma_T}\int_{\delta_T}^\infty 
\Bigl( 1 + v^2 \Bigr)^{1/4}
\exp \left( \frac{T\varphi_T}{8} v^2 \Bigl( 1 + v^2
\Bigr)^{-3/4} \right)\,dv.
\label{STARTDT}
\end{eqnarray}
Let $g$ and $h$ be the two functions defined on $\R^+$ by
$$
g(v) = \frac{v^2}{(1 + v^2)^{3/4}} \hspace{1cm} \text{and} \hspace{1cm}
h(v) = \frac{v^{3/2}}{(1 + v^2)^{3/4}}.
$$
One can observe that  $g$ and $h$ are both increasing functions on $\R^+$.
Moreover, as soon as $ v>\delta_T$, $g(v)=\sqrt{v} h(v)> \sqrt{v}h(\delta_T)$. In addition, for all
$v \in \R^+$, we also have
$$2^{3/4} g(v) \geq \min \left( v^2, \sqrt{v} \right). $$
Therefore, we obtain from \eqref{STARTDT} that
\begin{equation} 
\int_{|u| > s_T} |\Phi_T(u) |^2 \, du \leq 
\frac{ 8 C_{\ell}  T}{\gamma_T}
\exp\Bigl(\frac{T\varphi_T}{16}g(\delta_T)\Bigr) 
\int_{\delta_T}^\infty 
2^{1/4}\max(1, \sqrt{v})
\exp ( e_T \sqrt{v} )\,dv
\label{TEMPDT}
\end{equation}
where
$$
e_T=\frac{T\varphi_T}{16}h(\delta_T).
$$
The fact that $\varphi_T < 0$ leads to
\begin{eqnarray*}
\frac{T\varphi_T}{8}g(\delta_T)& \leq & \frac{T\varphi_T}{16}2^{3/4} g(\delta_T), \\
&\leq &  \frac{T\varphi_T}{16}\,\min \left( \delta_T^2,\sqrt{\delta_T} \right),\\
&\leq  & \frac{T\varphi_T}{16} \,\min \left(  \frac{4 c^2 }{\varphi_T^4} \frac{s_T^2}{\beta_T^2 }, 
\sqrt{ \frac{2 |c| }{ \varphi_T^2}} \sqrt{\frac{s_T}{|\beta_T|}} \right),\\
&\leq & \max \left( \frac{c^2 }{4 \varphi_T^3}, - \frac{\sqrt{2 |c| }}{ 16} \right)\,
\min \left( T\frac{s_T^2}{\beta_T^2 }, T \sqrt{\frac{s_T}{|\beta_T|}} \right), \\
& \leq &  - \mu C T^{\nu}
\end{eqnarray*}
where the positive constant $\mu$ in the last inequality is due to the boundeness of the terms in the max
and the power $T^{\nu}$ follows from assumption \eqref{HYPHYP}. Furthermore, for $T$ large enough,
the integral in \eqref{TEMPDT} is bounded by $1$. As a matter of fact, we have via straightforward calculation
on the Gamma function that
\begin{equation*} 
\int_{0}^\infty \max(1, \sqrt{v})
\exp ( e_T \sqrt{v} )\,dv \leq \frac{1}{e_T} \max\Bigl(1, -\frac{2}{e_T} \Bigr).
\end{equation*}
It is not hard to see from assumption \eqref{HYPHYP} that $e_T$ goes to $-\infty$ as $T$ tends to infinity,
which clearly implies that this integral is as small as one wishes. Then, we infer from \eqref{TEMPDT} 
that for $T$ large enough
\begin{equation} 
\int_{|u| > s_T} |\Phi_T(u) |^2 \, du \leq 
\frac{ 8 C_{\ell}  T}{\gamma_T}
\exp(-\mu C T^{\nu}) \leq \frac{ 8 C_{\ell}  T |\beta_T|\varphi^2_T}{2|c|}
\exp(-\mu C T^{\nu}).
\label{DDT}
\end{equation}
Finally, we deduce from \eqref{MAJDTCS}, \eqref{INTDT} and \eqref{DDT} that 
for $T$ large enough
\begin{equation*}
|D_T|^2 \leq \frac{ 8 C_{\ell}  T \varphi^2_T}{8 \pi |\alpha_T c|}
\exp(-\mu C T^{\nu})
\end{equation*}
which clearly implies that, for two positive constants $d$ and $D$,
$$
|D_T| \leq d ~  T~ \exp(-D T^{\nu})
$$
and completes the proof of Lemma \ref{L-GENE}. \demend

\subsection{Proof of Lemma \ref{L-BT1}.}
\begin{lem} \label{L-FC1}
For $c < - \theta$, the distribution of $U_T$ under $\p_T$ converges, as $T$ goes to infinity, to an
$\cN(0,1)$ distribution which means that $\Phi_T$ converges to $\Phi$ given by
$$
\Phi(u) = \exp\Bigl(-\frac{u^2}{2}\Bigr).
$$
Moreover, for any $p >0$, there exist integers $q(p)$, $r(p)$ and a sequence $(\varphi_{k,l})$ independent of $p$, such that, for $T$ large enough
\begin{equation} \label{DEVPHIC1}
\Phi_T(u)= \Phi(u)
\left[1+\frac{1}{\sqrt{T}}\sum_{k=0}^{2p}\sum_{l=k+1}^{q(p)}
\frac{\varphi_{k,l} u^l}{(\sqrt{T})^k}
+{\cal O}\Bigl(\frac{\max(1,|u|^{r(p)})}{T^{p+1}}\Bigr)\right]
\end{equation}
where $\sigma_c^2$ is given by \eqref{DEFACE1} and the remainder ${\cal O}$ is uniform as soon as $|u|\leq s T^{1/6}$ with $s>0$.
\end{lem}

\noindent{\bf Proof of Lemma \ref{L-FC1}.}
It follows from \eqref{E-DECFOND} that for all $k \in \N$, 
\begin{equation}
\label{DERLT} 
\cL_T^{(k)}(a_c) = \cL^{(k)}(a_c) + \frac{1}{T} \cH^{(k)}(a_c) + \frac{1}{T} \cR_T^{(k)}(a_c).
\end{equation}
Moreover, it is rather easy to see that for all $k \in \N$,\\ 
\begin{equation} \label{REST1}
\cR_T^{(k)}(a_c)={\cal O}(T^k\exp(2Tc)).
\end{equation}
One can observe that $\cL^{(1)}(a_c)=0$ and $\cL^{(2)}(a_c) =  \sigma_c^2 $
with $\sigma_c^2 $ given by \eqref{DEFACE1}. In addition, taking  $\beta_T = \sigma_c \sqrt{T}$, we also have
$$
T \,\left(\frac{iu}{\beta_T}\right)^{2}\!  \frac{\cL^{(2)}(a_c)}{2} = - \frac{u^2}{2}.
$$
Hence,  by a Taylor expansion, we find from \eqref{PHITOLT}, \eqref{DERLT} and \eqref{REST1} that for any $p >0$
\begin{equation*}
\log \Phi_T(u) = -\frac{u^2}{2} + T\sum_{k=3}^{2p+3}\left(\frac{iu}{\sigma_c\sqrt{T}}\right)^{k}\!  \frac{\cL^{(k)}(a_c)}{k!} + \sum_{k=1}^{2p+1} \left(\frac{iu}{\sigma_c \sqrt{T}}\right)^{k} \!  \frac{\cH^{(k)}(a_c)}{k!} + {\cal O}\Bigl(\frac{\max(1,u^{2p+4})}{T^{p+1}}\Bigr).
\end{equation*}
Finally, we deduce the asymptotic expansion \eqref{DEVPHIC1} by taking the exponential on both sides, remarking that, as soon as $|u|\leq s T^{1/6}$ with $s>0$, the quantity $u^l/(\sqrt{T})^k$ remains bounded in \eqref{DEVPHIC1}.
\demend

\noindent{\bf Proof of Lemma \ref{L-BT1}.}
In order to achieve the proof of the Lemma \ref{L-BT1}, let $s_T=s T^{1/6}$ with $s>0$ and $\beta_T= \sigma_c \sqrt{T}$.
As 
$$
\min \left( \frac{T s_T^2}{\beta_T^2}, \frac{T \sqrt{s_T}}{\sqrt{|\beta_T|}} \right) \geq C T^{1/3},
$$
the assumption \eqref{HYPHYP} of Lemma \ref{L-GENE} is clearly satisfied. Consequently, there exist two positive constants $d$ and $D$ such that
\begin{equation*}
|D_T|\leq d\exp(-D T^{1/3}).
\end{equation*}
Finally, we obtain \eqref{DEVBT1} from \eqref{CT} and \eqref{DEVPHIC1} together with standard calculations on the $\cN(0,1)$ distribution.
\demend

\subsection{Proof of Lemma \ref{L-BT2}.}

\begin{lem} \label{L-FC2}
For $c > \theta$, the distribution of $U_T$ under $\p_T$ converges, as $T$ goes to infinity, to the distribution
of $\gamma (N^2-1)$, where $N$ is an $\cN(0,1)$ random variable and
$$
\gamma=\cL^{\prime}(2(c-\theta))=\frac{(3c-\theta)}{2(\theta -2c)}
$$
which means that $\Phi_T$ converges to $\Phi$ given by
$$
\Phi(u) = \frac{\exp(-i\gamma u)}{\sqrt{1 - 2i\gamma u}}.
$$
Moreover, for any $p >0$, there exist integers $q(p)$, $r(p)$, $s(p)$ and a sequence $(\varphi_{k,l,m})$ independent of $p$, such that, for $T$ large enough
\begin{equation} 
\label{DEVPHIC2}
\Phi_T(u)= \Phi(u)\exp\left(-\frac{\sigma_c^2 u^2}{2T}\right) 
\left[1+\!\sum_{k=1}^p\sum_{l=k+1}^{q(p)}\sum_{m=0}^{r(p)}\frac{\varphi_{k,l,m}u^l}{T^k(1-2i \gamma u)^m}
+{\cal O}\Bigl(\frac{\max(1,|u|^{s(p)})}{T^{p+1}}\Bigr)\right]
\end{equation}
where $\sigma_c^2$ is given by \eqref{DEFACE2} and the remainder ${\cal O}$ is uniform as soon as $|u|\leq s T^{2/3}$ with $s>0$.
\end{lem}

\noindent{\bf Proof of Lemma \ref{L-FC2}.} It follows from \eqref{E-DECFOND} 	and \eqref{PHITOLT} that
\begin{equation} 
\label{LOC3}
\Phi_T(u) = \exp \left( \!T \!\left( \cL \left( a_T +  \frac{i u}{\beta_T} \right)\!-\! \cL(a_T) \right) + \cH \left( a_T +  \frac{i u}{\beta_T} \right)\!-\! \cH(a_T) 
+ \cR \left( a_T +  \frac{i u}{\beta_T} \right)\!-\! \cR(a_T) \right) 
\end{equation}
where $\beta_T=-T$. We shall focus our attention on each term of \eqref{LOC3}. First,
by virtue of Lemma \ref{lem-holo-r}, the term involving the remainder $\cR$ does not contribute
to the asymptotic expansion of $\Phi_T$. Next, by the very definition \eqref{E-DEFL} of $\cL$, the first term of \eqref{LOC3}
can be rewritten as
\begin{equation*}
T\left(\cL\Bigl(a_T+\frac{iu}{\beta_T}\Bigr)-\cL(a_T) \right)=-\frac{T}{2}\left( \frac{iu}{\beta_T}  - \varphi_T \left(\left( 1+\frac{iub_T}{\beta_T}\right)^{1/2}-1\right)\right)
\end{equation*}
where $\varphi_T=-\sqrt{\theta^2 +2 a_T c}$ and $b_T=2c/\varphi^2_T$. Consequently, as $b_T / \beta_T$ tends to 0,  
we have for all $p\geq 2$,
\begin{equation*}
\exp \left(  T\left(\cL\left(a_T+\frac{iu}{\beta_T}\right)-\cL(a_T) \right)\right)=\exp\left(
-\frac{iuT}{2\beta_T}  + \frac{T\varphi_T}{2} \sum_{k=1}^pl_k\Bigl(\frac{iub_T}{\beta_T}\Bigr)^k
+{\cal O}\Bigl(\frac{|u|^{p+1}}{T^{p+1}}\Bigr)\right)
\end{equation*}
where $l_k=(-1)^{k-1}(2k)!/((2k-1)(2^k k!)^2)$ which leads to
\begin{align}
 \exp \left( T  \left(\cL(a_T+\frac{iu}{\beta_T})-\cL(a_T) \right)\right) & \notag \\
= \exp\left( - iuc_T  - \frac{d_T u^2}{2T} \right) & \exp\left( \frac{T\varphi_T}{2} 
\sum_{k=3}^pl_k\Bigl(\frac{iub_T}{\beta_T}\Bigr)^k
+{\cal O}\Bigl(\frac{|u|^{p+1}}{T^{p+1}}\Bigr)\right)
\label{DEVLHARD}
\end{align}
where 
$$c_T=\frac{c-\varphi_T}{2\varphi_T} \hspace{1cm} \text{and}\hspace{1cm} 
d_T=-\frac{\varphi_T b_T^2}{8}. $$ 
For the second term of \eqref{LOC3}, we also have
by the very definition \eqref{E-DEFH} of $\cH$ 
\begin{equation*}
 \exp \left(  \cH\Bigl(a_T+\frac{iu}{\beta_T}\Bigr)-\cH(a_T) \right) = 
\left( \frac{\varphi_T + a_T +\theta}
{\varphi_T + \left( a_T + iu\beta_T^{-1} + \theta \right) \left( 1 + iu b_T\beta_T^{-1} \right)^{-1/2}}\right)^{1/2}\!\!.
\end{equation*}
Hence, we have for all $p\geq 2$,
\begin{align}
\exp \left(  \cH\Bigl(a_T+\frac{iu}{\beta_T}\Bigr)-\cH(a_T) \right) &  \notag \\
 = \frac{1}{\sqrt{f_T(u)}} &  \left(1 + g_T(u) u^2 +  h_T(u) \left( \sum_{k=2}^p h_k \Bigl(\frac{iub_T}{-T}\Bigr)^k
+{\cal O}\Bigl(\frac{|u|^{p+1}}{T^{p+1}}\Bigr) \right) \right)^{-1/2} 
\label{DEVHHARD}
\end{align}
where $h_k=(2k)!/(2^k k!)^2$ and 
\begin{eqnarray*}
f_T(u)&=& 1 - \frac{iu}{e_T} + \frac{(a_T+\theta)iub_T }{2e_T}, \\
g_T(u)&=& \frac{b_T }{2 T e_T f_T(u)}, \\
h_T(u)&=&  \frac{T(a_T+\theta) - iu}{e_Tf_T(u)},
\end{eqnarray*}
with $e_T=T(\varphi_T + a_T + \theta)$. One can easily check that, as $T$ goes to infinity, 
the limits of $b_T$, $c_T$, $d_T$, and $e_T$
are respectively given by $2c/(\theta - 2c)^2$, $\gamma$, $\sigma_c^2$, and $(\theta - c)/(3c - \theta)$
which implies that $f_T(u)$ converges to $1-2i \gamma u$.
Finally, we find via \eqref{DEVLHARD} and \eqref{DEVHHARD} the pointwise convergence
\begin{equation*}
\lim_{T \rightarrow \infty} \Phi_T(u)= \Phi(u)= \frac{\exp(-i \gamma u)}{\sqrt{1-2i \gamma u}}
\end{equation*}
while \eqref{DEVPHIC2} follows from the Taylor expansion of the exponential in \eqref{DEVLHARD} together 
with the Taylor expansion of the square root in \eqref{DEVHHARD}. 
\demend

\noindent{\bf Proof of Lemma \ref{L-BT2}.}
In order to complete the proof of the Lemma \ref{L-BT2}, let $s_T=s T^{2/3}$ with $s>0$ and $\beta_T= -T$.
It is not hard to see that
$$
\min \left( \frac{T s_T^2}{\beta_T^2}, \frac{T \sqrt{s_T}}{\sqrt{|\beta_T|}} \right) \geq C T^{1/3},
$$
which means that the assumption \eqref{HYPHYP} of Lemma \ref{L-GENE} is satisfied. 
Therefore, there exist two positive constants $d$ and $D$ such that
\begin{equation*}
|D_T|\leq d~T~\exp(-D T^{1/3}).
\end{equation*}
Finally, we deduce \eqref{DEVBT2} from \eqref{CT} and \eqref{DEVPHIC2} via a careful use of the contour integral lemma
for the Gamma function given in Lemma 7.3 of \cite{BR01}.
\demend

\subsection{Proof of Lemma \ref{L-BT3}.}

\begin{lem} \label{L-FC3}
For $|c| < \theta$ with $c\neq 0$, the distribution of $U_T$ under $\p_T$ converges, as $T$ goes to infinity, to the distribution
of $\gamma (N^2-1)$, where $N$ is an $\cN(0,1)$ random variable and
$$
\gamma=-\cL^{\prime}(0)=\frac{(\theta+c)}{2\theta}
$$
which means that $\Phi_T$ converges to $\Phi$ given by
$$
\Phi(u) = \frac{\exp(-i\gamma u)}{\sqrt{1 - 2i\gamma u}}.
$$
Moreover, for any $p >0$, there exist integers $q(p)$, $r(p)$, $s(p)$ and a sequence $(\varphi_{k,l,m})$ independent of $p$, such that, for $T$ large enough
\begin{equation} 
\label{DEVPHIC3}
\Phi_T(u)= \Phi(u)\exp\left(-\frac{\sigma_c^2 u^2}{2T}\right) 
\left[1+\!\sum_{k=1}^p\sum_{l=k+1}^{q(p)}\sum_{m=0}^{r(p)}\frac{\varphi_{k,l,m}u^l}{T^k(1-2i \gamma u)^m}
+{\cal O}\Bigl(\frac{\max(1,|u|^{s(p)})}{T^{p+1}}\Bigr)\right]
\end{equation}
where $\sigma_c^2$ is given by \eqref{DEFACE3} and the remainder ${\cal O}$ is uniform as soon as $|u|\leq s T^{2/3}$ with $s>0$.
\end{lem}

\noindent{\bf Proof of Lemma \ref{L-FC3}.} The proof is left to the reader inasmuch as it follows essentially the same lines
as those in the proof of Lemma \ref{L-FC2}. \demend

\noindent{\bf Proof of Lemma \ref{L-BT3}.} The proof of Lemma \ref{L-BT3} follows exactly the same
arguments as those in the proof of Lemma \ref{L-BT2}. The only notable thing to mention is that
we have to take in account twice the asymptotic behavior of $a_T$ because $a_T$ goes to zero as $T$ tends to infinity
and $a_T$ is also in the denominator of $C_T$. \demend

\subsection{Proof of Lemma \ref{L-BT4}.}

\begin{lem} \label{L-FC4}
For $c = -\theta$, the distribution of $U_T$ under $\p_T$ converges, as $T$ goes to infinity, to the distribution
of $\sigma_{\theta} N + \gamma_{\theta}(M^2-1)$, where $\sigma^2_{\theta}$ is given by \eqref{DEFAC4},
$N$ and $M$ are two independent $\cN(0,1)$ random variables and
$$
\gamma_{\theta}  = \frac{1}{2 \sqrt{\theta}}
$$
which means that $\Phi_T$ converges to $\Phi$ given by
$$
\Phi(u) = \frac{\exp \left( -i \gamma_{\theta} u \right)}{\sqrt{1 - 2i\gamma_{\theta} u}} \exp \left(- \frac{ \sigma_{\theta}^2 u^2}{2} \right).
$$
Moreover, for any $p >0$, there exist integers $q(p)$, $r(p)$, $s(p)$ and a sequence $(\varphi_{k,l,m})$ independent of $p$, such that, for $T$ large enough
\begin{equation} 
\label{DEVPHIC4}
\Phi_T(u)= \Phi(u) 
\left[1+\!\frac{1}{\sqrt{T}}\sum_{k=0}^{2p}\sum_{l=k+1}^{q(p)}\sum_{m=0}^{r(p)}\frac{\varphi_{k,l,m}u^l}{(\sqrt{T})^k(1-2i \gamma_\theta u)^m}
+{\cal O}\Bigl(\frac{\max(1,|u|^{s(p)})}{T^{p+1}}\Bigr)\right]
\end{equation}
where the remainder ${\cal O}$ is uniform as soon as $|u|\leq s T^{1/6}$ with $s>0$.
\end{lem}

\noindent{\bf Proof of Lemma \ref{L-FC4}.} The proof is left to the reader inasmuch as it follows essentially the same lines
as those in the proof of Lemma \ref{L-FC2}. \demend

\noindent{\bf Proof of Lemma \ref{L-BT4}.} The proof of Lemma \ref{L-BT4} follows exactly the same
arguments as those in the proof of Lemma \ref{L-BT2} with a careful use of the contour integral lemma
for the Gamma function given in Lemma 7.3 of \cite{BR01}. \demend

\section{Appendix D: Technical Lemmas.}
\setcounter{equation}{0}

\subsection{Statement of the results.}

The effective domain $\Delta_{T,c}$ of the normalized cumulant generating function
$\cL_T$ was previously calculated in Lemma 4.1 of \cite{FP99}. It is an open interval 
such that $\widetilde{\Delta}_{T,c} \subset \Delta_{T,c} \subset \overline{\Delta}_{T,c}$ where
$$\widetilde{\Delta}_{T,c} = \Bigl\{ a \in \R,  ~\theta^2 +2ac > 0,~a + \theta < \sqrt{\theta^2+2ac}\coth(T\sqrt{\theta^2+2ac}) \Bigr\}$$
and
$$\overline{\Delta}_{T,c} = \Bigl\{ a \in \R,  ~\theta^2 +2ac+\frac{\pi^2}{T^2} > 0,~a + \theta < \sqrt{\theta^2+2ac}\coth(T\sqrt{\theta^2+2ac}) \Bigr\}.$$
Denote $\cD_{T,c}=\{z\in \C, \re(z) \in \Delta_{T,c}\}$ and $\cD_{c}=\{z\in \C, \re(z) \in \Delta_{c}\}$.

\begin{lem} \label{L-LTCA}
For $T$ large enough, $\Phi_T$ belongs to $\dL^2(\R)$. More precisely, for 
$T$ large enough and for any $(a,u) \in \R^2$ such that $(a+iu) \in \cD_{T,c}$,
\begin{equation}
\label{MAJPHIT}
\Bigl|\exp \left( T \left( {\mathcal L}_T (a+iu) -{\mathcal L}_T(a) \right) \right) \Bigr|^2 \leq 4\ell(a,c,\theta)  
\Bigl(1+ \frac{4c^2u^2}{\varphi^4(a)}\Bigr)^{1/4} \exp\left(T\frac{c^2 u^2}{2\varphi^3(a)} 
\Bigl( 1 +\frac{4c^2 u^2}{\varphi^4(a)} \Bigr)^{-3/4}\right)
 \end{equation}
where
\begin{equation} 
\ell(a,c,\theta) 
= \max\left(1, \frac{|\varphi(a) + \theta |}{|\varphi(a) |} \right)
\max\left(1, \frac{|\varphi(a) + 2c-\theta |}{|\varphi(a) |} \right).
\label{DEFLMAJ}
\end{equation}
\end{lem}

\subsection{Proof of Lemma \ref{L-LTCA}.}

The key point is to make use of a complex counterpart of the main  decomposition \eqref{E-DECFOND}, which means
\begin{equation} \label{E-DECFONDCOMP}
\cL_T(z)=\cL(z) + \frac{1}{T} \cH(z) + \frac{1}{T} \cR_T(z)
\end{equation}
where $\cL$, $\cH$ and $\cR_T$ are respectively given by \eqref{E-DEFL}, \eqref{E-DEFH} and \eqref{E-DEFR}. 
In order to make these expressions meaningful, we have to take care about the definitions. We shall denote the principal determination of the logarithm defined on $\C \setminus ]-\infty, 0]$ by
\begin{align*}
\loga{z} = \log |z| + i \Arg(z),
\end{align*}
where
\begin{eqnarray*}
\Arg(z) = \left \{ \begin{array}{lll}
\arcsin \left( \frac{\im(z)}{|z|} \right) & \text{ if }   \re(z) \geq 0, \vspace{1ex} \\
   \arccos \left( \frac{\re(z)}{|z|} \right) & \text{ if }  \re(z) <0,~\im(z)>0, \vspace{1ex} \\
    \!\! -\arccos \left( \frac{\re(z)}{|z|} \right)    & \text{ if } \re(z) <0,~\im(z)<0. 
\end{array}  \right.
\end{eqnarray*}
We also introduce the analytic function defined for all $z \in \C$ with $Re(z)>0$, by
\begin{equation*}
\sqrt{1+ z}= \sqrt{|1+z|} \exp \Bigl( \frac{i}{2} \, \Arg(1+z) \Bigr).
\end{equation*}
It is not hard to see that
\begin{equation} \label{RELATION}
\re(\sqrt{1+z})= \frac{1}{\sqrt{2}} \sqrt{ | 1+z | + 1+\re(z)}.
\end{equation}
The proof of Lemma \ref{L-LTCA} follows from the conjunction of three lemmas, each one involving
the functions $\cL$, $\cH$ and $\cR_T$.

\begin{lem}\label{L-HOLOL}
The function $\cL$ given, for all $z \in \C$, by
$$
\cL(z)= -\frac{1}{2}\Bigl(z+ \theta -\varphi(z)\Bigr) \hspace{1cm} \text{where} \hspace{1cm}
\varphi(z)=- \sqrt{\theta^2 + 2zc}
$$
is differentiable on the domain $\cD_c$. Moreover, for all $a \in \Delta_c$ and $u \in \R$, we have
\begin{equation}
\label{MAJFFL}
\Bigl| \exp \left(T( \cL(a+iu)-\cL(a)) \right)\Bigr|^2 \leq \exp \left( T\frac{c^2 u^2}{4\varphi^3(a)} 
\Bigl( 1 +\frac{4c^2 u^2}{\varphi^4(a)} \Bigr)^{-3/4}\right).
\end{equation}
\end{lem}

\noindent{\bf Proof of Lemma \ref{L-HOLOL}.}
For all $z \in \C$ such that $\re(z) \in \Delta_c$, $\varphi(z)$ is well defined. Hence, $\varphi$ is differentiable on $\cD_c$ and the same
is also true for $\cL$. In addition, we have
\begin{equation*}
\cL(a+iu)-\cL(a) = -\frac{1}{2} \left( iu - \varphi(a+iu) + \varphi(a)\right)
\end{equation*}
which clearly implies that
\begin{equation}
\Bigl| \exp \left(T( \cL(a+iu)-\cL(a)) \right)\Bigr| \leq 
\exp \left(\frac{T}{2}( \re ( \varphi(a+iu) - \varphi(a))) \right). 
\label{MAJFFL1}
\end{equation}
Moreover, we also have
\begin{equation*}
\varphi(a+iu) - \varphi(a) = \varphi(a) \left( \sqrt{ 1 + \frac{2icu}{\varphi^2(a)} } - 1 \right).
\end{equation*}
We deduce from \eqref{RELATION} with $z = 2icu/\varphi^2(a)$, that
$$
\re \left( \varphi(a+iu) - \varphi(a) \right)
=
\frac{\varphi(a)}{\sqrt{2}} \left( \sqrt{\sqrt{1 + \frac{4c^2u^2}{\varphi^4(a)}} + 1 } - \sqrt{2} \right).
$$
Keeping in mind that $\varphi(a) < 0$, we infer from the elementary inequality
$$
\sqrt{\sqrt{1 + x} + 1 } - \sqrt{2} \geq \frac{x}{ 4 \sqrt{2} ( 1 + x )^{3/4}}
$$
which is true as soon as $x \geq 0$, that
$$
\re \left( \varphi(a+iu) - \varphi(a) \right)
\leq 
\frac{c^2 u^2}{2\varphi^3(a)} \left( 1 +\frac{4c^2 u^2}{\varphi^4(a)} \right)^{-3/4}.
$$
Finally, it ensures via \eqref{MAJFFL1} that for all $a \in \Delta_c$ and $u \in {\mathbb R}$,
\begin{equation*}
\Bigl| \exp \left(T( \cL(a+iu)-\cL(a)) \right)\Bigr|^2 \leq \exp \left( T\frac{c^2 u^2}{4\varphi^3(a)} 
\Bigl( 1 +\frac{4c^2 u^2}{\varphi^4(a)} \Bigr)^{-3/4}\right).
\end{equation*}
which ends the proof of Lemma \ref{L-HOLOL}. \demend

\begin{lem}\label{L-HOLOH}
The function $\cH$ given, for all $z \in \C$, by
$$
\cH(z)= - \frac{1}{2}\log \left( \frac{1}{2} ( 1 +h(z) ) \right)
 \hspace{1cm} \text{where} \hspace{1cm}
h(z)= \frac{(z+\theta)}{\varphi(z)}
$$
is differentiable on the domain $\cD_c$. Moreover, for all $a \in \Delta_c$ and $u \in \R$, we have
\begin{equation}
\label{MAJFFH}
\Bigl| \exp (\cH(a+iu) - \cH(a)) \Bigr|^2 \leq \ell(a,c,\theta)  \left(1+ \frac{4c^2u^2}{\varphi^4(a)}\right)^{1/4}.
\end{equation}
\end{lem}

\noindent{\bf Proof of Lemma \ref{L-HOLOH}.}
First of all, it follows from \eqref{E-DEFH} that
\begin{equation}
\label{MAJFFH1}
\Bigl| \exp \left(\cH(a+iu)-\cH(a) \right)\Bigr|^2 = \left| \frac{1+h(a)}{1+h(a+iu)} \right|. 
\end{equation}
We claim that for $z \in {\mathbb C}$ such that $\re(z) \in \Delta_c$, $1+ h(z) \in {\mathbb C}\setminus ]-\infty, 0]$.
Assume by contradiction that this is not true, which means that one can find $\lambda \in [0,+\infty[$ such that
$$
1+h(z)=-\lambda.
$$
Since, $\varphi^2(z)=\theta^2 +2cz$, $\varphi(z)$ is a root of the quadratic equation
$$
\varphi^2(z)+2c(1+\lambda) \varphi(z) - \theta^2 +2c\theta=0.
$$
Its discriminant is non-negative as
$$
4(c-\theta)^2 + 4c^2 \lambda^2 + 8 c^2 \lambda \geq 0.
$$
One can observe that $c$, $\theta$ and $\lambda$ are reals numbers which implies that $\varphi(z)$ is also a real number as well as $z$. 
Consequently, $z$ belongs to $\Delta_c$ and $1+ h(z)> 0$ which contradicts the assumption. 
This allows us to say that $\cH$ is differentiable on $\cD_c$. We are now in position to prove 
inequality \eqref{MAJFFH}.
Since $\varphi^2(z)=\theta^2 +2cz,$ for $z\in {\mathbb C}$ such that $\re(z) \in \Delta_c$, we have
\begin{equation}\label{DECOH}
1+h(z)= \frac{(\varphi(z) + \theta)(\varphi(z) +2c - \theta)}{2c\varphi(z)}.
\end{equation}
For all $z \in {\mathbb C}$ such that $\re(z) \in \Delta_c,$ and for all $\alpha \in {\mathbb R}$, we clearly have  
$$
\left|\varphi(z)+\alpha\right|^2= |\varphi(z)|^2 +\alpha^2 + 2\alpha \re(\varphi(z)).
$$
Assume that $a$ belongs to $\Delta_c$ and let $u \in \R$. We already saw that
\begin{equation*} 
\varphi(a+iu)  = \varphi(a) \sqrt{1+ \frac{2ic u}{\varphi^2(a)}}
\end{equation*}
which leads to
\begin{equation}
\label{MAJFFH2} 
|\varphi(a+iu)|^2 = \varphi^2(a) \left| 1+ \frac{2icu}{\varphi^2(a)} \right|=
\varphi^2(a) \left(1+ \frac{4c^2u^2}{\varphi^4(a)}\right)^{1/4}.
\end{equation}
On the other hand, it also follows from \eqref{RELATION} that
\begin{equation}
\label{REAL}
\re(\varphi(a+iu)) = \frac{\varphi(a)}{\sqrt{2}}~\sqrt{ 1+\left |1+\frac{2icu}{\varphi^2(a)}\right| }=
\frac{\varphi(a)}{\sqrt{2}}~\sqrt{1+\sqrt{ 1 + \frac{4c^2 u^2}{\varphi^4(a)}}}.
\end{equation}
Consequently,
$
|\varphi(a+iu)|^2 = 2(\re(\varphi(a+iu)))^2-\varphi^2(a)
$
which implies that
\begin{equation} \label{MODULE1}
\Bigl|\varphi(a+iu) + \alpha \Bigr|^2=\varphi^2(a) \left(\frac{2(\re(\varphi(a+iu)))^2}{\varphi^2(a)}-1 +\frac{\alpha^2}{\varphi^2(a)} +
\frac{2 \alpha}{\varphi(a)}~\frac{\re(\varphi(a+iu))}{\varphi(a)}\right).
\end{equation}
By introducing the function
\begin{equation*}
g_{\alpha}(x)
= \left( x + \frac{\alpha}{\sqrt{2}} \right)^2 + \frac{\alpha^2}{2} -1,
\end{equation*}
we deduce from \eqref{REAL} together \eqref{MODULE1} with $\beta=\alpha /\varphi(a)$, that
\begin{equation}
\label{MODULE2}
\Bigl|\varphi(a+iu) + \alpha \Bigr|^2= \varphi^2(a)~g_{\beta}\left( \sqrt{1+\sqrt{ 1 + \frac{4c^2 u^2}{\varphi^4(a)}}}\right).
\end{equation}
Furthermore, one can easily check from straighforward calculations that for all $x \geq \sqrt{2}$,
\begin{eqnarray}
\label{MING}
g_{\beta}(x) \geq 
\left \{ \begin{array}{ll}
(\beta +1)^2 & \text{ if }   \beta \in [-2,0], \vspace{1ex} \\
1    & \text{ otherwise.}
\end{array}  \right.
\end{eqnarray}
Therefore, we infer from \eqref{MODULE2} and \eqref{MING} that for all
$a\in \Delta_c$ and for all $u \in \R$,
\begin{eqnarray*}
|\varphi(a+iu) + \alpha \Bigr | \geq 
\left \{ \begin{array}{ll}
|\varphi(a) | |\beta +1| & \text{ if }   \beta \in [-2,0], \vspace{1ex} \\
|\varphi(a) |    & \text{ otherwise }
\end{array}  \right.
\end{eqnarray*}
which clearly implies that
\begin{equation}
\label{YANG}
\frac{1}{|\varphi(a+iu) + \alpha |} \leq \frac{1}{ |\varphi(a) |} \max\left(1, \frac{|\varphi(a) |}{|\varphi(a) + \alpha |} \right). 
\end{equation}
We shall make use of inequality \eqref{YANG} with $\alpha=\theta$ and $\alpha=2c-\theta$. One can observe that, as long as
$a\in \Delta_c$, the value of $\varphi(a) + \alpha \neq 0$. Finally, it follows from the conjunction of
\eqref{MAJFFH1}, \eqref{DECOH}, \eqref{MAJFFH2}, and \eqref{YANG} that for all
$a\in \Delta_c$ and for all $u \in \R$,
\begin{eqnarray*}
\Bigl| \exp \left(\cH(a+iu)-\cH(a) \right)\Bigr|^2 &\leq &
\frac{|\varphi(a+iu)|}{|\varphi(a)|} \max\left(1, \frac{|\varphi(a) + \theta |}{|\varphi(a) |} \right)
\max\left(1, \frac{|\varphi(a) + 2c-\theta |}{|\varphi(a) |} \right), \\
&\leq &  \ell(a,c,\theta)  \left(1+ \frac{4c^2u^2}{\varphi^4(a)}\right)^{1/4}
\end{eqnarray*}
which completes the proof of Lemma \ref{L-HOLOH}. \demend

\begin{lem}\label{L-HOLOR}
For $T$ large enough, the function $\cR_T$ given, for all $z \in \C$, by
$$
\cR_T(z) = -\frac{1}{2}  \log \left( 1+\frac{1-h(z)}{1 +h(z)} \exp(2 \varphi(z) T) \right)
$$
is differentiable on the domain $\cD_{T,c}$. Moreover, for all 
$(a,u) \in \R^2$ such that $a + iu \in \cD_{T,c}$, we have
\begin{equation}
\label{MAJFFR}
\Bigl| \exp (\cR_T(a+iu) - \cR_T(a)) \Bigr|^2 \leq 4.
\end{equation}
\end{lem}

\noindent{\bf Proof of Lemma \ref{L-HOLOR}.}
First of all, we deduce from \eqref{E-DEFR} that
\begin{equation}
\label{MAJFFR1}
\Bigl| \exp \left(\cR_T(a+iu)-\cR_T(a) \right)\Bigr|^2 = \left| \frac{1+r(a)\exp(2\varphi(a) T) }{1+r(a+iu)\exp(2\varphi(a+iu) T) } \right|
\end{equation}
where the function $r$ given, for all $z \in \C$, by
$$
r(z) = \frac{1-h(z)}{1 +h(z)}.
$$
We already saw from \eqref{DECOH} that
\begin{equation*}
1+h(z)= \frac{(\varphi(z) + \theta)(\varphi(z) +2c - \theta)}{2c\varphi(z)}.
\end{equation*}
Hence,
\begin{equation*}
1-h(z)= \frac{(\theta - \varphi(z)  )(\varphi(z) -2c + \theta)}{2c\varphi(z)}
\end{equation*}
which implies that
\begin{equation}
\label{DECOR}
| r(z) | = \left| \frac{\varphi(z)-\theta}{\varphi(z)+\theta} \right|  \left| \frac{\varphi (z) -2c + \theta }{\varphi(z) + 2c -\theta} \right|.
\end{equation}
Moreover, for all $(a,u) \in \R^2$ such that $a + iu \in \cD_{T,c}$, 
\begin{equation*}
\Bigl|r(a+iu)\exp(2\varphi(a+iu) T)\Bigr|^2=|r(a+iu)|^2 \exp(4 \re(\varphi(a+iu)) T).
\end{equation*}
We recall from \eqref{REAL} that
\begin{equation*}
\re(\varphi(a+iu)) = 
\frac{\varphi(a)}{\sqrt{2}}~\sqrt{1+\sqrt{ 1 + \frac{4c^2 u^2}{\varphi^4(a)}}}.
\end{equation*}
The key point here is that $\re(\varphi(a+iu))$ is always negative. Consequently, via the same lines as in the proof of Lemma \ref{L-HOLOH},
we obtain that for $T$ large enough and for all $(a,u) \in \R^2$ such that $a + iu \in \cD_{T,c}$, 
\begin{equation}
\label{MAJFFR2}
\Bigl|r(a+iu)\exp(2\varphi(a+iu) T)\Bigr| \leq \frac{1}{2}.
\end{equation}
Finally, \eqref{MAJFFR} follows from \eqref{MAJFFR1} and \eqref{MAJFFR2}. \demend

\noindent{\bf Proof of Lemma \ref{L-LTCA}.} Lemma \ref{L-LTCA} immediately follows from \eqref{E-DECFONDCOMP}
together with the conjunction of Lemmas \ref{L-HOLOL}, \ref{L-HOLOH} and \ref{L-HOLOR}. \demend

\ \vspace{-1ex} \\
{\bf Acknowledgements.}                              
The first author is deeply grateful to Alain Rouault for fruitful discussions on
a preliminary version of the manuscript.

\bibliographystyle{plain}
\bibliography{BCSExplo}

\end{document}